\documentclass[prx,twocolumn,amsmath,amssymb,notitlepage]{revtex4-1}

\usepackage{url}
\usepackage{hyperref}
\hypersetup{
    colorlinks,
    citecolor=blue,
    filecolor=black,
    linkcolor=blue,
    urlcolor=black
}

\usepackage{mathrsfs}
\usepackage{verbatim}
\usepackage{graphicx}	
\usepackage{dcolumn}	
\usepackage{bm}		
\usepackage{bbm,upgreek}
\usepackage{mathtools}
\usepackage{algpseudocode}

\usepackage{color}
\usepackage[cal=boondoxo]{mathalfa}
\usepackage{palatino}

\usepackage{dsfont}
\usepackage{cancel}
\usepackage{xcolor}
\usepackage{harpoon}
\usepackage{accents}
\usepackage{xr}
\usepackage[margin=0.78in,footskip=0.25in]{geometry}

\newcommand{\expect}[1]{\underset{#1}{\mathbb{E}}}

\newcommand*{\Scale}[2][4]{\scalebox{#1}{$#2$}}%

\newcommand*{\inliner}[1]{{\smash{\Scale[0.85]{#1}}}}

\definecolor{Gray}{gray}{0.85}

\usepackage{xcolor,colortbl}

\definecolor{maroon}{cmyk}{0,0.87,0.68,0.32}

\DeclareFontFamily{U}{euc}{}
\DeclareFontShape{U}{euc}{m}{n}{<-6>eurm5<6-8>eurm7<8->eurm10}{}
\DeclareSymbolFont{AMSc}{U}{euc}{m}{n}
\DeclareMathSymbol{\umu}{\mathord}{AMSc}{"16}

\DeclareSymbolFont{AMSb}{U}{msb}{m}{n}
\DeclareMathSymbol{\N}{\mathbin}{AMSb}{"4E}
\DeclareMathSymbol{\Z}{\mathbin}{AMSb}{"5A}
\DeclareMathSymbol{\R}{\mathbin}{AMSb}{"52}
\DeclareMathSymbol{\Q}{\mathbin}{AMSb}{"51}
\DeclareMathSymbol{\I}{\mathbin}{AMSb}{"49}

\begin{document}

\title{On the correspondence between thermodynamics and inference}

\author{Colin H. LaMont and Paul A. Wiggins}
\affiliation{Departments of Physics, Bioengineering and Microbiology, University of Washington, Box 351560.\\ 3910 15th Avenue Northeast, Seattle, WA 98195, USA}

\email{pwiggins@uw.edu}\homepage{http://mtshasta.phys.washington.edu/}

\begin{abstract}
We  expand upon a natural analogy between Bayesian statistics and statistical physics in which sample size corresponds to inverse temperature.
This analogy motivates the definition of two novel statistical quantities: a learning capacity and a Gibbs entropy.
The analysis of the learning capacity, corresponding to the heat capacity in thermal physics, leads to new insight into the mechanism of learning and explains why some  models  have anomalously-high learning performance.
We explore the properties of the learning capacity in a number of examples, including a sloppy model.
Next, we  propose that the Gibbs entropy provides a natural device for counting distinguishable distributions in the context of Bayesian inference. We use this device to define a \textit{generalized principle of indifference} (GPI) in which every distinguishable model is assigned equal \textit{a priori} probability. This principle results in a new solution to a long-standing problem in Bayesian inference: the definition of an objective or uninformative prior.
A key characteristic of this new approach is that it can be applied to analyses where the model dimension is unknown and circumvents the automatic rejection of higher-dimensional models in Bayesian inference.
\end{abstract}


\maketitle

\section{Introduction}
%



Despite an intensifying interest in applications of machine learning to the analysis of big data, fundamental questions remain about the mechanism of learning and the development of efficient learning algorithms.
In this paper, we explore the phenomenology of learning by exploiting a correspondence between Bayesian inference and statistical mechanics.
This correspondence has been previously described by Jaynes, Balasubramanian, and many others \cite{Jaynes1957information,Jaynes2003,barndorff1997statistics,Balasubramanian1997,nishimori2001statistical,mezard2009information} and  methods from statistical physics have been adapted to statistical calculations \cite{doucet2000sequential,hoffman2014no,chen2014stochastic,minka2001expectation,xing2002generalized,minka2005divergence,niinimaki2013annealed,liang2001real,inoue2001application,mezard2006reconstruction,zdeborova2016statistical}.
Motivated by the success of this previous work, we propose to exploit the correspondence at a more conceptual level. By using the canonical bridge between statistical mechanics and thermodynamics, we define statistical analogues to the standard thermodynamic potentials and properties of a system. We then explore the statistical properties of these new analogues.
The correspondence identifies two novel statistical quantities, a \textit{learning capacity} and the \textit{Gibbs entropy} which give new physical insight into the mechanism of learning and defines a novel Bayesian learning algorithm, respectively.

\begin{figure*}
  \centering

    \includegraphics[width=.95\textwidth]{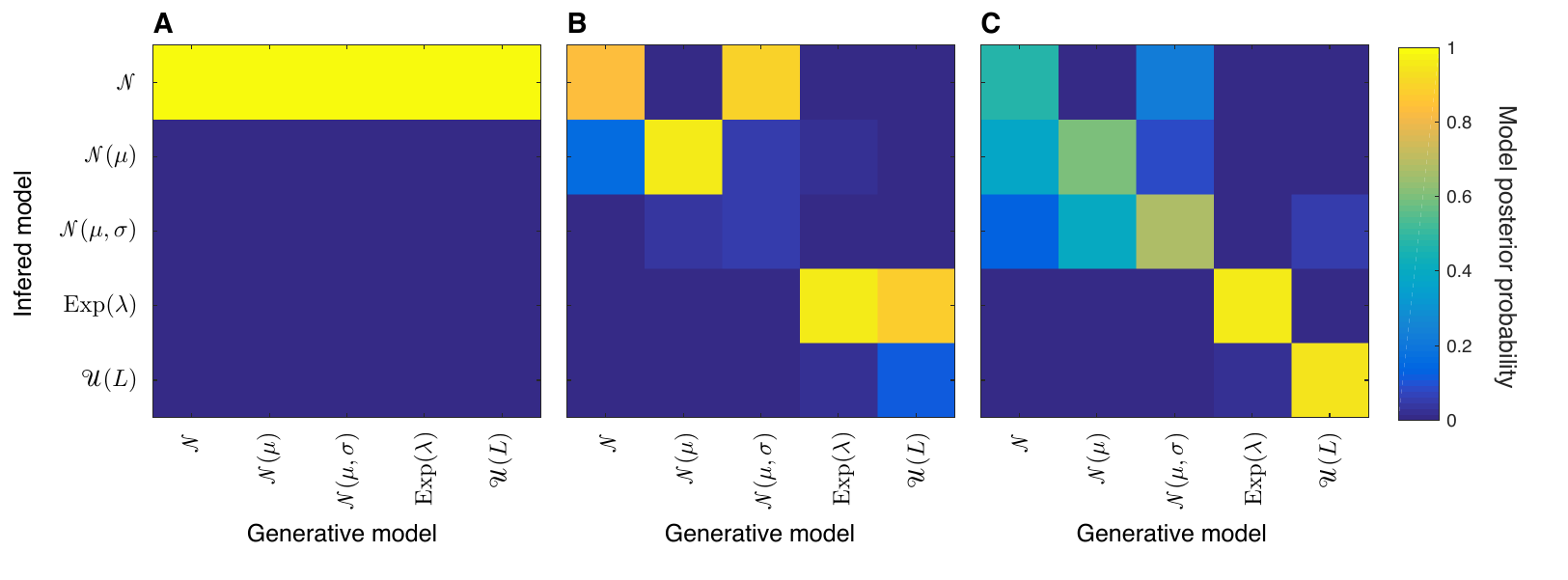}
      \caption{ \textbf{Bayesian inference on model identity.} The posterior of model identity ($y$ axis) was computed for datasets generated by each model ($x$ axis).  \textbf{Panel A: Objective prior.} The non-compactness of the parameter manifolds implies automatic rejection of all the higher-dimensional models. Since the model ${\cal N}$ is parameter-free, it has posterior probability of 1 for all datasets, regardless of the fit. \textbf{Panel B: Revised objective prior.} To avoid this undesirable anomaly, we tune the prior parameter support to result in a reasonable posterior model probability. (See Tab.~\ref{tab:models}.) Inference is no longer objective, as the posterior probabilities depend on how this tuning is performed. One representative  plot of posterior probabilities to shown. In general, inference cannot be expect to identify the generative model unless the KL divergence is so large as to make the prior irrelevant.     \textbf{Panel C: GPI prior.} Using the GPI prior, the generative model had the highest posterior probability as expected. When ${\cal N}$ is the generative distribution, the data is realizable in both ${\cal N}(\mu)$ and ${\cal N}(\mu,\sigma)$ as well. The fact that these model have lower posterior probability reveals that there remains a natural mechanism for the selection of parsimonious models using the  GPI prior.
      \label{Fig:inference}}
\end{figure*}

\subsection{Model complexity}
How does model complexity affect learning and why do some models have anomalously good learning performance? 
These are questions of great topical interest due to the increasing application of machine learning algorithms and the analysis of extremely complex models in the context of systems biology and other fields (\textit{e.g.}~\cite{Machta:2013hl}).
 The analysis of a novel \textit{learning capacity} (corresponding to the heat capacity) reveals a natural connection between the equipartition theorem in statistical mechanics \cite{pathria1996statistical}, which  predicts the thermal energy as a function of the \textit{number of degrees of freedom of a physical system}, and the information loss in prediction as a function of the \textit{number of degrees of freedom of a model}. This connection provides physical insight into why there are universal properties of learning systems that are  independent of the detailed functional dependence of the likelihood functions or the learning algorithm.  
 
In spite of these universal limits, it has long been known that some high-dimensional models learn anomalously well. These models have been termed \textit{sloppy}  \cite{Machta:2013hl}. 
We demonstrate that the learning capacity both provides a natural definition for the sloppiness phenomenon and identifies a mechanism  for the anomalously-high predictive performance, a statistical analogue of the well-known freeze-out mechanism of statistical mechanics.
We speculate that this mechanism is responsible for the anomalously-high predictive performance of   high-dimensional models more generally.

\subsection{Prior selection}

The correspondence also provides new insight into another important question: How do you define an objective or uninformative prior? Prior selection has been a subject of debate in Bayesian inference from its inception to the current day \cite{Kass1996}.
A key application of the correspondence is to translate insights  into improved learning algorithms. In this light, we propose that the Gibbs entropy provides a natural device for determining model multiplicity, \textit{i.e.}~counting indistinguishable distributions in the context of statistical inference.
This interpretation allows us to define a \textit{generalized principle of indifference} (GPI) for selecting a prior in the absence of \textit{a priori} information on the parameters or models. The GPI unifies a number of known, but seemingly unconnected objective Bayesian methods \cite{Gelman2014}, while also providing an algorithm applicable  in contexts where existing approaches fail \cite{Berger2012}. 

\subsection{The Lindley-Bartlett Paradox}

\label{sec:lindbartintro}

The most important advantage of the GPI approach over existing approaches is in the context of models with unknown dimension (\textit{i.e.}~model selection \cite{akaike1773,BurnhamBook,Gelman2014}). Here, the use of  uninformative priors can often result in the automatic rejection of higher-dimensional models, a phenomenon called the Lindley-Bartlett paradox \cite{Bartlett1957,Lindley1957}. To demonstrate this phenomenon, we  analyze simulated datasets where the generative distribution is known and test the performance of inference. As we will demonstrate, the canonical Bayesian approach to inference leads to difficulties in generating a meaningful posterior probability on model identity.

\medskip
\noindent
\textit{Numerical experiment:} We consider five competing models: three realizations of the normal model (known mean and variance, unknown mean and known variance, and unknown mean and variance), the exponential model with unknown rate and the uniform distribution model with unknown end-point. 
We  generate datasets from all fives models.  We then perform inference on the model identity for each dataset using an uninformative prior.  
A detailed definition of the likelihood functions, priors, parameter support and generative model parameters are provided in the Appendix (\ref{Sec:CanBayes}) and Tab.~\ref{tab:models}.


If the canonical Bayesian approach is interpreted literally,  the posterior probability of the parameter-free normal model is one, regardless of which distribution was used to generate the data! (See Fig.~\ref{Fig:inference}A.)  This  approach ensures that the smallest model is always selected, however much data is accumulated and however poorly this model fits the data. Although this scenario would appear paradoxical, it perfectly logical from a Bayesian perspective:  Due to the non-compact manifold, a measure zero fraction  of the parameter space is consistent with the data. The relative merits of Bayesian and Frequentist inference in this context were first debated by Lindley and Bartlett in the 1950s \cite{Bartlett1957,Lindley1957}. 


A number of \textit{ad hoc} modifications to this naive Bayesian procedure are possible to avoid the automatic rejection of larger models. These include some formal methods: \textit{e.g.}~the use of empirical or variational Bayes, pseudo-Bayes factors, partitioning of the data to train the prior \textit{etc.}
Many practical-minded Bayesians will revise their prior beliefs, presumably remembering \textit{a posteri} that they knew \textit{a priori}.
 We use this approach in Fig.~\ref{Fig:inference}B. 
In this context not only does inference fail to identify the generative distribution, but it is also depends sensitively on \textit{ad hoc} decisions made in the analysis, which we would argue is not acceptable in the context of objective scientific analysis.

How can priors be defined over models to give sensible and robust results? The use of the proposed GPI prior circumvents the paradox by naturally assigning a mutually-consistent weighting, not only between parameter values, but between model families, irrespective of model dimension. (The GPI approach is shown in Fig.~\ref{Fig:inference}C.)
Our approach is asymptotically consistent with a powerful but non-Bayesian approach to inference on model identity: the use of the pseudo-Bayes factor \cite{Gelman2014,gelfand1994bayesian,Spiegelhalter2002,Vehtari2012,BURNHAM2004}. Thus the GPI approach provides a novel solution to one of the oldest problems in statistics: the specification of an objective prior.

\subsection{Outline}

The paper is organized as follows: In Sec.~\ref{sec:corr}-\ref{sec:analogues}, we define the correspondence and compute the thermodynamic potentials and properties of inference in the analytically tractable large-sample-size limit and use these results to deduce the statistical meaning of each quantity. In Sec.~\ref{sec:learningCap}, we explore the statistical properties of the  learning capacity. In Sec.~\ref{sec:Gibbs}, we use the Gibbs entropy to define a generalized principle of indifference and an objective  prior (the GPI prior). In Sec.~\ref{sec:appGPI}, we compute the GPI prior in a number of examples.  In Sec.~\ref{sec:lind}, we discuss how the GPI prior circumvents the Lindley-Bartlett paradox.

\begin{table*}[t]
\setlength{\tabcolsep}{1em} 
{\renewcommand{\arraystretch}{1.5}
\begin{tabular}{|ll|c|ll|}
\hline
{\bf Thermodynamics} 	&    &   & {\bf Statistics}  &  \\
{\bf Quantity:}        &  {\bf Interpretation:} & & {\bf Quantity:}     &  {\bf Interpretation:} \\              
\hline
\hline
$\beta = { T}^{-1}$ & Inverse temperature & $\leftrightarrow$ &
	
			$N$ &  Sample size \\
			${\bm \theta}$ & State variables/vector  & $\leftrightarrow$& ${\bm \theta}$ & Model parameters \\
			$X^N$ & Quenched disorder  & $\leftrightarrow$ & $X^N$ & Observations \\ 
			${ E}_X({\bm \theta})$ & State energy  & $\leftrightarrow$ & $\hat{H}_X({\bm \theta})$ & Cross entropy estiamtor \\ 
			${ E}_0$ & Disorder-averaged ground state energy  & $\leftrightarrow$ & $H_0$ & Shannon entropy \\ 
		    $\rho({\bm \theta})$ & Density of states  & $\leftrightarrow$ & $\varpi({\bm \theta})$ & Prior  \\
 ${Z}$  & Partition function  & $\leftrightarrow$    & $Z$    & Evidence  \\     
 		    ${Z}^{-1} \rho\ \exp -\beta { E}_X$ & Normalized Boltzmann weight  & $\leftrightarrow$ & $\varpi({\bm \theta}|X^N)$ & Posterior  \\	
		    \hline	  
${ F} = -\beta^{-1}\log {Z}$	& Free energy & $\leftrightarrow$& $F = -N^{-1} \log Z$ 
        	& Minus-log-evidence  \\
$ { U} =   \partial_\beta \beta { F}$ & Average energy  & $\leftrightarrow$&   $U = \partial_N N F$ 	& Minus-log-prediction  \\ 
 ${ C} = -\beta^2 \partial^2_\beta \beta { F}$ & Heat capacity & $\leftrightarrow$ & ${ C} = -N^2 \partial^2_N N {F}$ 	& \textit{Learning capacity} \\
        
	 ${ S} = \beta^2\partial_\beta F$ & Gibbs entropy & $\leftrightarrow$      	& ${S} = N^2\partial_N F$ & \textit{Statistical Gibbs entropy} \\
            \hline
\end{tabular} 
\caption{{\bf Thermodynamic-Bayesian correspondence.} The top half of the table lists the correspondences that can be determined directly from the definition of the marginal likelihood as the partition function. The lower half of the table lists the implied thermodynamic expressions and their existing or proposed statistical interpretation. \label{Table:stat}}}
\end{table*}

\section{Methods}

\subsection{Defining the correspondence} 
\label{sec:corr}

We assume that a true parameter value ${\bm \theta}_0$ is drawn from a known prior distribution $\varpi({\bm \theta})$. We observe $N$ samples $x^N\equiv\{x_1,...,x_N\}$ which are distributed like  $q(x|{\bm \theta}_0)$:
\begin{equation}
X_i\sim q(\cdot|{\bm \theta}_0),
\end{equation}
where we use capital $X$ to denote random variables and the symbol $\sim$ to denote \textit{distributed like}.
For simplicity, we will assume that the observations are independent and identically distributed, but the approach can be generalized.


The correspondence between statistical physics and Bayesian inference is clearest when expressions are written in terms of the empirical estimator of the cross-entropy:
\begin{eqnarray}
\hat{H}({\bm \theta}) &\equiv& -\left<\log q(X|{\bm \theta})\right>_{X\sim x^N},\\
&\equiv& -N^{-1} \sum_{i=1}^N \log q(x_i|{\bm \theta}), \label{eqn:emp2}
\end{eqnarray}
where the angle brackets represent an expectation over the variable in the subscript, which in this case is an empirical expectation over the observed data $\inliner{x^N}$ (Eq.~\ref{eqn:emp2}). 
The  marginal likelihood (\textit{i.e.}~evidence) can be written \cite{Balasubramanian1997}:
\begin{align}
Z(x^N)&\equiv \int_{\bm \Theta}\!\! {\rm d}{\bm \theta}\ \varpi({\bm \theta})\,e^{-N\hat{H}},\label{eqn:evidence}
\end{align}
which can  be directly compared to the partition function in the canonical ensemble\footnote{In the canonical ensemble, the probability of the occupancy of a microstate $i$ is $p_i = Z^{-1} g_ie^{-\beta E_i}$ where $Z$ is the partition function, $E_i$ is the energy, $g_i$ is the degeneracy and $\beta^{-1} \equiv k_BT$.} \cite{Jaynes1957information,Jaynes2003,barndorff1997statistics,Balasubramanian1997,nishimori2001statistical,mezard2009information}. The model parameters ${\bm \theta}$ correspond to  the variables that define the physical state vector, the cross entropy estimator $\inliner{\hat{H}({\bm \theta})}$ corresponds to the energy $ E({\bm \theta})$, the prior $\varpi({\bm \theta})$ corresponds to  the density of states $\rho({\bm \theta})$. 
The data $x^N$ is  quenched disorder\footnote{Quenched disorder refers to variables in the system that are randomly generated when the system is assembled but remain constant in time. These variables remain fixed in expectations over thermal fluctuations.}
 in the physical system \cite{Balasubramanian1997}.
The sample size $N$ is identified with the inverse temperature $\beta\equiv { (k_BT)}^{-1}$ \cite{Balasubramanian1997}. (Henceforth, we will set $k_B\equiv 1$.) This assignment is natural in the following sense: At small sample size $N$, many parameter values are consistent with the data, in analogy with the large range of states $\bm \theta$ occupied at high temperature.
In contrast, at large sample size $N$ the parameter values consistent with the data are tightly localized around the true value, in analogy to a statistical system at low temperature where  only states $\bm \theta$ in very close proximity to the ground state are occupied.
We note that choosing $T^{-1} \leftrightarrow N$ is only one of at least two proposals for the identification of the temperature. See the Appendix (\ref{App:alttemp}).

\subsection{Application of thermodynamic identities}
To extend the previously proposed correspondence, we follow the standard prescriptions from statistical mechanics to compute thermodynamic potentials, properties, and variables for the system  \cite{gibbs1902elementary,pathria1996statistical}. These are shown in the lower half of Tab.~\ref{Table:stat}.
The thermodynamic quantities depend on the particular realization of the data $X^N$.
In the current context we are interested in the expectation over this \textit{quenched disorder} (\textit{i.e.}~data).
We define the disorder average with an overbar:
\begin{align}
\overline{f}(N) \equiv \left<\right. f(X^N,{\bm \theta}_0) \left.\right>_{X,{\bm \theta}_0},
\end{align}
where $X\sim q(\cdot|{\bm \theta}_0)$ and ${\bm \theta}_0\sim \varpi$. 

%

\label{sec:res}

\section{Results}

\subsection{Models}
Our immediate interest in the next few sections is not performing inference but rather exploring the properties of the statistical counterparts of well-understood thermodynamic quantities. The similarity between the statistical and thermodynamic quantities suggests that these novel statistical quantities may have an analogous interpretation to their thermodynamic counterparts. We will motivate this hypothesis by comparing the properties of statistical models in the the large-sample-size limit to the properties of a free particle (\textit{i.e.}~a gas).

\subsubsection{Free particle} 

In the context of statistical mechanics (and thermodynamics), we will analyze a  classical free particles in $K$ dimensions in  the canonical ensemble.  The particle has internal energy $E_0$ and a phase-space density of states $h^{-K}$, where $h$ is the Planck constant. The particle is confined to a $K$-cube with side-length $L$. We define a critical temperature $T_0$  at which the de Broglie wavelength is equal to $L$: 
\begin{equation}
T_0 \equiv \beta_0^{-1} \equiv h^2/2\pi mL^2.
\end{equation}
The thermodynamics quantities are straightforward to compute and are summarized in Tab.~\ref{tab:results}. The Free energy is:
\begin{equation}
F =  { E}_0 + \textstyle \frac{K}{2\beta } \log \frac{\beta}{\beta}_0, \label{eqn:FPF}
\end{equation}
which is written in terms of the critical inverse temperature $\beta_0$, a parameter which holds all the information about both the  density of states as well as the geometry of the system.

\begin{table*}
\setlength{\tabcolsep}{1.6em} 
{\renewcommand{\arraystretch}{1.5}
\begin{tabular}{|lc|l|l|l|}
\hline
\textbf{} & & \textbf{K-Dimensional} & \textbf{K-Dimensional} & \textbf{K-Dimensional}  \\
\textbf{} & & \textbf{Free Particle} & \textbf{Regular Model} & \textbf{Singular Model}  \\
\hline
\hline
Free energy & $\overline{F}$  & ${ E}_0 + \textstyle \frac{K}{2\beta } \log \frac{\beta}{\beta}_0$   & $H_0 + {\textstyle \frac{K}{2N}} \log{\textstyle \frac{N}{N_0}}\ \ +{\cal O}(N^{-1})$ & $H_0 + \textstyle \frac{\gamma}{2N} \log N\ \ \, +{\cal O}(N^{-1})$  \\
Average energy & $\overline{U}$   & ${ E}_0 + \textstyle \frac{K}{2\beta }$   & $H_0 + \textstyle \frac{K}{2N }\ \ \ \ \ \ \ \ \ \ \ \ \ \ +{\cal O}(N^{-2})$ &  $H_0+\textstyle \frac{\gamma}{2N}\ \ \ \ \ \ \ \ \, \ \ \ \ +{\cal O}(N^{-2})$  \\
Heat capacity & $\overline{C}$  &  $\textstyle \frac{K}{2 }$  &   $\textstyle \frac{K}{2 }\ \ \ \ \ \ \ \ \ \ \ \ \ \ \ \ \ \ \ \ \ \ \ \ \ \, +{\cal O}(N^{-1})$  & $\textstyle \frac{\gamma}{2}\ \ \ \ \ \ \ \ \ \ \ \ \ \ \ \ \ \ \ \ \ \ \ \ \ \, +{\cal O}(N^{-1})$  \\
Gibbs entropy & $\overline{S}$  &  ${\textstyle \frac{K}{2}}\left(1-\log \frac{\beta}{\beta_0}\right)$  & ${\textstyle \frac{K}{2}}\left(1-\log \frac{N}{N_0}\right)\ \, +{\cal O}(N^{-1})$   & $-{\textstyle \frac{\gamma}{2}}\log N\ \ \ \ \ \ \ \ \ \ \ \ +{\cal O}(N^{0})$  \\
\hline
\end{tabular} 
\caption{{\bf Thermodynamic-Bayesian correspondence.} The thermodynamic quantities of a $K$-dimensional free particle with ground-state energy $E_0$ are compared to a $K$-dimensional regular model.
Inspection reveals that the $\beta/N$ dependence of free particle is identical to a regular model to order $N^{-1}$.
For the singular model, the learning coefficient is $\gamma \le K$. The special case of  $\gamma = K$ is a regular model. \label{tab:results}}}
\end{table*}

\subsubsection{Large-sample-size limit in regular and singular models}

In the context of statistics, we will compare and contrast two classes of models: \textit{regular} and \textit{singular}. 
Models are called \textit{singular} when parameters are \textit{structurally unidentifiable}, defined as the absence of a one-to-one map between the space of candidate distribution functions ($q$) and the parameter manifold. In other words there exists some parameter ${\bm \theta}_1$ such that
\begin{equation}
q(x|{\bm \theta}_1)=q(x|{\bm \theta}_2)\ \ \ \ \ {\rm for}\ \ \ \ \   {\bm \theta}_1\ne {\bm \theta}_2.
\label{eqn:unident}
\end{equation}
A model is singular when the  unidentifiability cannot be removed by re-parameterization \cite{watanabe2009}. In this case, the Fisher Information Matrix (Eq.~\ref{fisherMat}) contains at least one zero eigenvalue at ${\bm \theta}_1$ \cite{watanabe2009}. In contrast, in a \textit{regular model}, all parameters are identifiable, the parameter manifold is continuous, and the Fisher information matrix  is therefore positive definite in a suitable parameter coordinate system. 

\subsection{Definition of thermodynamic analogues}

There are known  asymptotic results for the evidence (partition function) in large-sample-size limit of both regular and singular models on continuous parameter manifolds \cite{watanabe2009}. It is therefore straightforward to compute the thermodynamic quantities in this limit once $Z$ is known. These results are shown in Tab.~\ref{tab:results}. 
A comparisons between the properties of a free particle and a regular model reveal an identical structure to leading order in $N$ (or $\beta$).



\label{sec:analogues}

\subsubsection{Free energy} The model that maximizes the evidence and therefore minimizes $F$ is selected in the canonical approach to Bayesian model selection \cite{Gelman2014}. 
Since, a relation between the partition function and Bayesian evidence has long been discussed \cite{Jaynes1957information,Jaynes2003,barndorff1997statistics,Balasubramanian1997,nishimori2001statistical,mezard2009information}, the definition of $F$ is not particularly novel.

For a regular model, $F$ breaks up into two parts:
\begin{equation}
 \overline{F} = H_0 + {\textstyle \frac{K}{2N}} \log{\textstyle \frac{N}{N_0}} +{\cal O}(N^{-1}),\label{eqn:bayesF}
\end{equation}
to order $\inliner{N^{-1}}$ where the dependence in the prior is absorbed into a \textit{critical sample-size} $N_0$. We call $N_0$ the \textit{critical sample size} because $N=N_0$ corresponds to the sample size at which the thermodynamic properties of inference change, as we will shall discuss shortly (Sec.~\ref{sec:normmod1}). Due to the dependence of $F$ on the critical sample size $N_0$, the evidence is clearly dependent on the choice of prior, if only logarithmically.


From the perspective of statistical mechanics, a direct comparison between this expression  (Eq.~\ref{eqn:bayesF}) and the free energy of a free particle (Eq.~\ref{eqn:FPF}) allows the reader intuitively understand the motivation for the correspondence defined in Sec.~\ref{sec:corr}: $H_0$ corresponds to the ground-state energy and the second term in Eq.~\ref{eqn:bayesF} is an entropic contribution to the free energy.

\subsubsection{Average energy}
The thermodynamic prescription for computing the average energy involves a derivative with respect to temperature (Tab.~\ref{Table:stat}):
\begin{equation}
U \equiv \partial_NNF.
\end{equation}
In the context of a discrete temperature, we will formally interpret this derivative using a finite difference definition:
\begin{equation}
\partial_Nf(N) \rightarrow f(N)-f(N-1).
\end{equation}
Such finite-difference approximations are already implicit to statistical mechanics, where we take derivatives with respect to many variables  which are in fact discrete (\textit{e.g.}~particle number, energy, \textit{etc.}). In the context of statistics where there are $N$  independent choices for reducing the sample size by one sample, it is convenient to define the finite difference derivative by averaging over the choices:
\begin{equation}
U(x^N) \equiv -\left<\log q(x_i|x^{\ne i})\right>_{i=1..N},\label{eqn:loocv}
\end{equation}
where  $\inliner{q(x_i|x^{\ne i})\equiv Z(x^N)/Z(x^{\ne i})}$ is known as the posterior-predictive distribution \cite{Gelman2014}. The RHS is a well-known statistical object: the Leave-One-Out-Cross-Validation (LOOCV) estimator of model performance. See the Appendix (\ref{App:crossval}).
The statistical interpretation of average energy $U$ is therefore the minus-expected-predictive-performance of the model (\textit{e.g.}~\cite{hastie01statisticallearning}).

To explore the correspondence to the free particle, we compute the averaged energy for the regular model. $U$ can be written as the sum of two contributions: 
\begin{equation}
\overline{U} = H_0 + \textstyle \frac{K}{2N }+{\cal O}(N^{-2}),
\end{equation}
to order $\inliner{N^{-2}}$.
The first term  $H_0$  corresponds to a ground-state energy. The second term  corresponds to the thermal energy in a physical system. From a statistical perspective, the term represents the information loss associated with predicting a new observation $X$ using parameters estimated from the training set $\inliner{x^{N}}$ rather than the true parameter ${\bm \theta}_0$. 
This predictive information loss is often called the \textit{Generalization Error}, defined (\textit{e.g.}~\cite{hastie01statisticallearning}):
\begin{equation}
{\rm GE} \equiv  H_0-\overline{U}.
\end{equation}
In a regular model, GE is \cite{watanabe2009}: 
\begin{equation}
{\rm GE} = -\textstyle \frac{K}{2N }+{\cal O}(N^{-2}),
\end{equation}
which depends only on the model dimension $K$ and sample size $N$ but it is independent of the detailed structure of the model (\textit{i.e.}~independent of the likelihood function $q$). 

This universal generalization error has a well-known thermodynamic analogue in the equipartition theorem: \textit{There is a half $k_BT$ of thermal energy per harmonic degree of freedom} \cite{pathria1996statistical}. In a statistical context, there is a universal generalization error of $\inliner{\frac{1}{2N}}$ per degree of freedom in the model. This universal property of the generalization error is known in statistics  (\textit{e.g.}~\cite{watanabe2009}) and can be interpreted as the mechanism by which the Akaike Information Criterion (AIC) estimates the predictive performance \cite{akaike1773,BurnhamBook}. However, the connection between this result and the equipartition theorem had not yet been described.



\subsubsection{Learning capacity}

To study the generalization error associated with learning from a finite-sized sample, it is natural to study the statistical quantity corresponding to the heat capacity.
The heat capacity measures the rate of increase in thermal energy with temperature ($\inliner{\overline{C}}$ in Tab.~\ref{Table:stat}).
The statistical analogue of the heat capacity, a \textit{learning capacity}:
\begin{equation}
C \equiv -N^2\partial_N U, 
\end{equation}
is a measure of the rate of increase in predictive performance with sample size. 

To explore the correspondence to the free particle, we compute the learning capacity for a regular model:
\begin{equation}
\overline{C} = \textstyle\frac{1}{2} K+ {\cal O}(N^{-1}), \label{eqn:equipart}
\end{equation}
as implied by the equipartition theorem. To learn how this analogy generalizes to a generic statistical model, we use the large-sample-size limit asymptotic expression for the Bayesian evidence for a singular model from Ref.~\onlinecite{watanabe2009} to compute the learning capacity. (See Tab.~\ref{tab:results}.)  
Like the normal model, the learning capacity for a singular model has the equipartition form but with an effective complexity:
\begin{equation}
\overline{C} = \textstyle\frac{1}{2} K_{\rm eff}+ {\cal O}(N^{-1}), \label{eqn:equipart}
\end{equation}
where $K_{\rm eff}=\gamma$ is the learning coefficient defined by Watanabe \cite{watanabe2009}. 
%
A regular model is a special case of this expression where $\gamma = K$, the dimension of the parameter manifold. The learning capacity is a novel statistical object defined by the correspondence. We will explore its properties in detail in Sec.~\ref{sec:learningCap}.


. 

\subsubsection{The Gibbs entropy} 
In physics, the Gibbs entropy generalizes the Boltzmann formula: $S = \log \Omega$ where $\Omega$ is the number of accessible states. We propose that the Gibbs entropy has the analogous meaning in the context of Bayesian statistics: The Gibbs entropy is the log-number of models consistent with the data. The correspondence implies that the statistical analogue to the Gibbs entropy is defined:
\begin{equation}
S(x^N) \equiv N(U-F),\label{eqn:loocv2}
\end{equation}
in analogy to statistical mechanics.


To explore the correspondence to the free particle, we compute the Gibbs entropy for a regular model:
\begin{equation}
\overline{S} = {\textstyle \frac{K}{2}}\left(1-\log \textstyle\frac{N}{N_0}\right)+{\cal O}(N^{0}).
\end{equation}
When the data is informative to the parameter values, the number of models consistent with the data is reduced as the sample size grows. As a result the Gibbs entropy is always negative for a normalized prior and becomes increasingly negative as the sample size $N$ grows. 

Like the learning capacity, the Gibbs entropy is a novel statistical object defined by the correspondence which we shall argue provides a natural mechanism for defining an objective prior in which all models are assigned equal \textit{a priori} weight. We will explore its properties in detail in Sec.~\ref{sec:Gibbs}.

\subsection{Examples of the Learning Capacity} 
\label{sec:learningCap}

\begin{figure}
  \centering
    \includegraphics[width=0.48\textwidth]{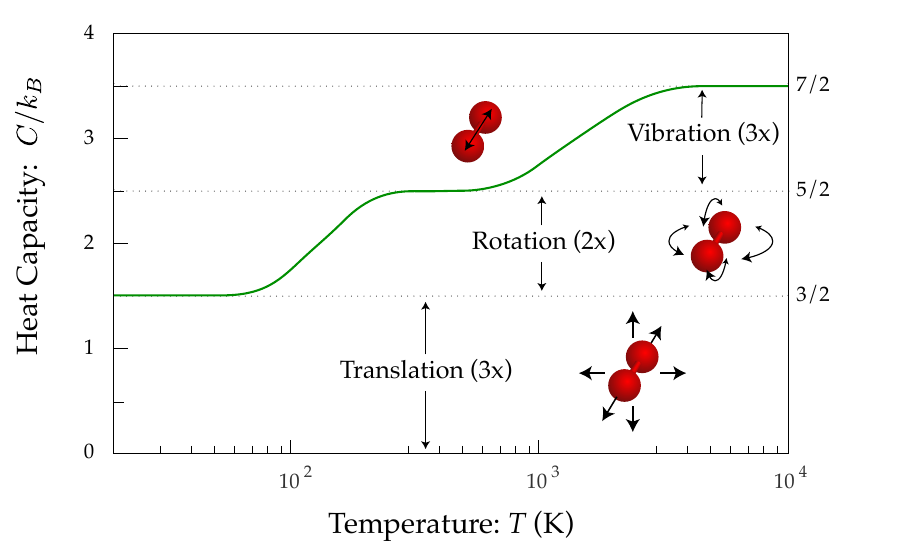}
      \caption{{\bf The failure of equipartition.} \textbf{Low-temperature freeze-out in a quantum system.} The heat capacity is plotted as a function of temperature. Equipartition predicts that the reduced heat capacity should be constant, equal to half the degrees of freedom in the system. Plateaus can clearly be observed at half-integer values, but the number of degrees of freedom is temperature dependent due to the discrete nature of quantum energy levels. At low temperature, some degrees of freedom are frozen out since the first excited state is thermally inaccessible.  This discrete topology of the energy levels implies anharmonicity in the potential and therefore failure of the equipartition theorem. 
      \label{Fig:DU1}}
\end{figure}

\label{sec:learningCap}

\begin{figure*}
  \centering
    \includegraphics[width=0.95\textwidth]{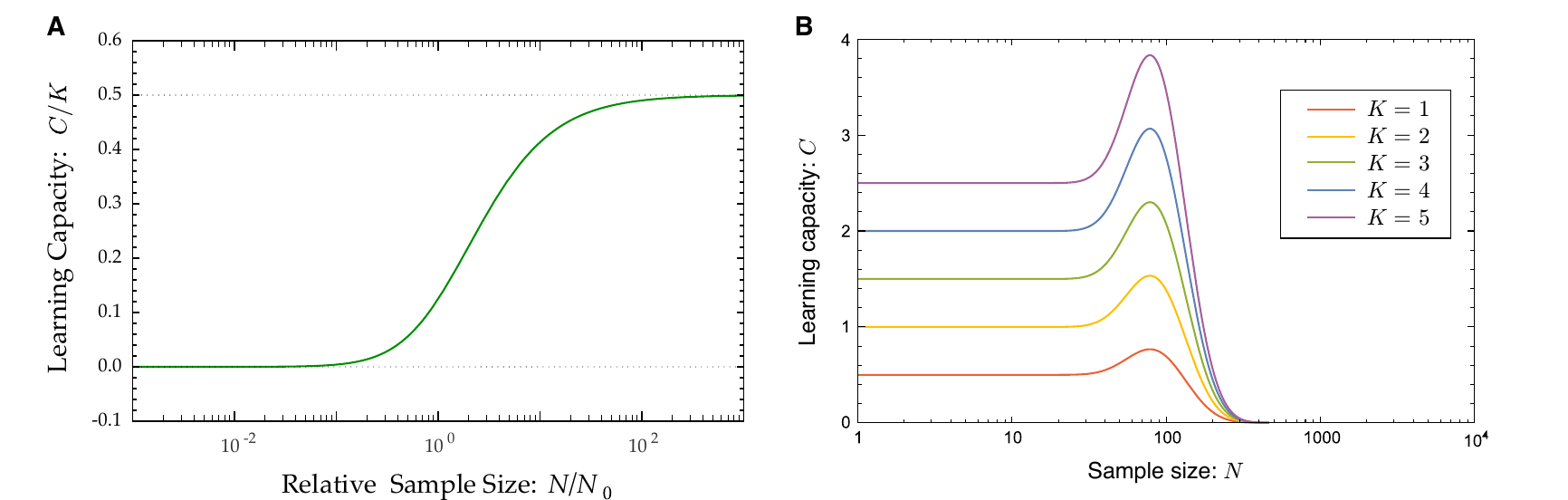}
      \caption{{\bf The failure of equipartition.} \textbf{Panel A: ``High-temperature freeze-out'' in the Learning Capacity.} Analogous to the statistical mechanics system, the statistical learning capacity transitions between half integer plateaus, reflecting a temperature-dependent number of degrees of freedom. At low sample size $N$ (high temperature), the parameters are completely specified  by model constraints (the prior) and therefore the parameters do not contribute to the learning capacity. At large sample size $N$, the parameters become data dominated and therefore the learning capacity is predicted by equipartition ($\frac{1}{2}K$). \textbf{Panel B: Low-temperature freeze out in the Learning Capacity.} The learning capacity of a normal model with an unknown D-dimensional mean $\inliner{\vec{\mu}\in {\mathbb Z}^D}$ and variance $\sigma^2 = 15$. For statistical uncertainty $\delta \mu \gg 1$, the learning capacity is predicted by equipartition since the discrete nature of the parameter manifold cannot be statistically resolved. For $\delta \mu \ll 1$, there is no statistical uncertainty in the parameter value (due to the discreteness of $\mu$) and the degrees of freedom freeze out, giving a learning capacity of zero.
      \label{Fig:DU2}}
\end{figure*}


In this section, we investigate the phenomenology of the learning capacity in a series of simple examples.  
We shall demonstrate that the learning capacity can show large deviation from the equipartition  limit.
From a physical perspective, this is no surprise. The failure of the equipartition theorem is a  well understood phenomenon in physics where degrees of freedom can become anharmonic  at both high and low temperature, altering their contribution to the heat capacity \cite{pathria1996statistical}. See Fig.~\ref{Fig:DU1}. We will analyze statistical models where analogous transitions occur as a function of sample size.

%

\subsubsection{ High-temperature freeze out}

It is well known in statistical physics that degrees of freedom can  become irrelevant at high temperature. For instance, the position degrees of freedom of a gas do not contribute to the heat capacity  \cite{Reif1967}. Exactly this type of phenomenon occurs in inference as well where a significant fraction of the degrees of freedom of a model can become anharmonic and are not data dominated. We present two simple examples of this fairly general phenomenon.

\subsubsection{ A normal model with an informative prior}
\label{sec:normmod1}
The canonical model in a statistical context is the normal model (a Gaussian distribution), many of whose properties generalize to more generic models in the large-sample-size limit.

\medskip
\noindent
\textit{Model:} We define a normal model on a $D$-dimensional observational space with \textit{unknown mean} and \textit{known variance} $\sigma^2$. The likelihood function is:
\begin{equation}
q(\vec{x}|{\bm \theta}) \equiv (2\pi \sigma^2)^{-D/2}\ \exp [-\textstyle\frac{1}{2\sigma^2}(\vec{x}-\vec{\mu}\,)^2], \label{eqn:normalmodel}
\end{equation}
with parameters ${\bm \theta} \equiv (\vec{\mu})$.
The true parameter ${\bm \theta}_0$ is drawn from a $K$-dimensional normal distribution:
\begin{equation}
\varpi({\bm \theta}) \equiv (2\pi \sigma_\varpi^2)^{-D/2}\ \exp [-\textstyle\frac{1}{2\sigma^2_\varpi}(\vec{\mu}-\vec{\mu}_\varpi\,)^2], \label{eqn:normalmodelprior}
\end{equation}
where $\mu_\varpi$ and $\sigma_\varpi$ are assumed to be known hyper-parameters.  
In close analogy to the inverse critical temperature of the free particle, it will be convenient to define a critical sample size:
\begin{equation}
N_0 \equiv \sigma^2/\sigma^2_\varpi, \label{eqn:critss}
\end{equation}
where $N_0$ can be understood as the number of previous observations $x^{N_0}$ required to determine the prior in the absence of other information.

\medskip
\noindent
\textit{Analysis:} The learning capacity can be computed analytically:  
\begin{eqnarray}
\overline{C} = \textstyle\frac{K}{2(1+N_0/N)^2},
\end{eqnarray}
as shown in the Appendix (\ref{sec:mmumkvp}).
At large sample size, the learning capacity is equal to the equipartition expression (Eq.~\ref{eqn:equipart}). At small sample size (high temperature), the prior determines the parametrization, $\inliner{\overline{C}\rightarrow 0}$ and the model appears anomalously predictive.  (See Figs.~\ref{Fig:DU1} and \ref{Fig:DU2}A.) 


\subsubsection{Exponential mixture models.}

\label{sec:expmix}

In the previous example, there are fairly transparent constraints applied to the parameters in the form of a prior which reduce the learning capacity. Our next example illustrates another \textit{higher-temperature} freeze-out phenomenon, but here the mechanism is model singularity: A zero mode appears in the Fisher information matrix.

\medskip
\noindent
\textit{Model:} We analyze the exponential mixture model which has previously been identified as \textit{sloppy model} by Transtrum, Machta and coworkers using a criterion defined by the distribution of the eigenvalues of the Fisher information matrix \cite{Machta:2013hl}. Consider a model for the lifetime of a mixed population consisting of several different chemical species $I$ with different transition rates. Both the transition rates ($k_I$) and the relative abundance of the species ($p_I$) are unknown. For an $m$ species model, the likelihood function for the lifetime $t$ is:
\begin{eqnarray}
q(t|{\bm \theta}) &\equiv& \sum_{I=1}^mp_I\, k_I\, {\rm e}^{-k_It},
\end{eqnarray}
with parameters:
\begin{eqnarray}
{\bm \theta} &\equiv& \left( \begin{array}{ccc} p_1 & ... & p_m \\ k_1 & ... & k_m \end{array} \right), 
\end{eqnarray}
subject to the constraint: $\inliner{\sum_I p_I = 1}$ and we apply improper prior $\varpi({\bm\theta}) = 1$. For simplicity, we analyze the smallest model with a singularity ($m=2$) to facilitate the numerical Bayesian marginalization. The exponential mixture model is singular since parameter $k_I$ is unidentifiable for $p_I=0$ and $p_1$ is unidentifiable for $k_1=k_2$. (See Eq.~\ref{eqn:unident}.)  

\medskip
\noindent
\textit{Analysis:} We compute the learning capacity at two locations in parameter manifold, at the singularity (${\bm \theta}_S$) and far from it (${\bm \theta}_R$):
\begin{equation}
{\bm \theta}_S = \left( \begin{array}{cc} 1 & 0 \\ 1 & 10 \end{array} \right) \ \ \ \ \ \ {\text{and}} \ \ \ \ \ {\bm \theta}_R = \left( \begin{array}{cc} {\textstyle\frac{1}{2}} & {\textstyle\frac{1}{2}} \\ 1 & 10 \end{array} \right).
\end{equation}
The learning capacity is computed numerically for $N=100$ observations with distribution $\inliner{T^N\sim q(\cdot|{\bm \theta})}$:
\begin{equation}
\overline{C}({\bm \theta}) = \begin{cases}
0.61, & {\bm \theta}= {\bm \theta}_S \\
1.5, & {\bm \theta}= {\bm \theta}_R
\end{cases}. 
\label{eqn:Llc}
\end{equation}
Far from the singularity (${\bm \theta}_R$), the equipartition theorem predicts the learning capacity ($\dim/2$) whereas close to the singularity (${\bm \theta}_S$), where the model is effectively described by only a single parameter ($k_1$), the learning capacity reflects this smaller effective model dimension.

\subsubsection{Low-temperature freeze out}

The freeze out phenomenon owes its name to the loss of harmonic degrees of freedom at low temperature in physical systems. Here the discrete nature of the quantum energy states plays an essential role in the physics. As the temperature becomes too low to populate the lowest-energy excited state, this degree of freedom no longer contributes to the heat capacity. Again, the analogous phenomenon is found in inference. In this context, the discrete nature of space is typically the result of a discrete parameter.

\label{sec:discrete}

\medskip
\noindent
\textit{Model:} To explore the low-temperature freeze-out phenomenon, consider a normal model with an unknown discrete mean. 
 The likelihood for the $D$-dimensional normal model is defined in Eq.~\ref{eqn:normalmodel}. As before, the parameters are ${\bm \theta} = (\vec{\mu})$, but with the mean constrained to have an integer values: $\vec{\mu}\in \mathbb{Z}^D$. 

\medskip
\noindent
\textit{Analysis:} The learning capacity can be computed analytically. See the Appendix (\ref{App:discra}). The learning capacity is plotted in Figs.~\ref{Fig:DU1} and \ref{Fig:DU2}B.

To discuss the phenomenology, it is useful to define a frequentist \textit{statistical resolution} with respect to parameter coordinate $\theta^i$:
\begin{equation}
\delta \theta^i(N) \equiv N^{-\frac{1}{2}}\sqrt{[\bm I^{-1}]^{ii}}, \label{eqn:resolution}
\end{equation}
in terms of the Fisher information matrix $\bm I$ (Eq.~\ref{fisherMat}), which is a naturally covariant symmetric tensor on the continuous parameter manifold \cite{Amari1985}. $\delta \theta^i(N)$ is the width of the posterior in the large-sample-size limit. For the normal model, $\delta \mu = \sigma/\sqrt{N}$. For a regular model with discrete parameters in the large sample size limit, the learning capacity is:
\begin{equation}
\overline{C} = \begin{cases}
\textstyle\frac{1}{2}K, & \Delta \theta^i \ll \delta \theta^i\ \text{for all }i \\
0, & \Delta \theta^i \gg \delta \theta^i\ \text{for all }i 
\end{cases} \label{eqn:lcdiscrete}
\end{equation}
where $\Delta \theta^i$ is the lattice spacing for parameter coordinate $\theta^i$. The physical interpretation is clear: At large sample size, the system condenses into a single state.  Therefore the corresponding degrees of freedom freeze out, and no longer contribute to the learning capacity. At small sample size, the discrete nature of the parameter manifold cannot be resolved, and the parameter manifold is effectively continuous. The learning capacity therefore assumes the equipartition value (provided that the sample size is large enough such that the information is effectively harmonic yet small enough so the discrete nature of the parameter manifold cannot be resolved).

\subsubsection{Learning capacity for additional models}

In the Appendix, we provide a series of other examples of learning capacity analyses to further explore its phenomenology.   
In Appendix Secs.~\ref{app:norm} and \ref{app:Exp}, we analyze a range of exactly tractable (but non-Gaussian) models. These analyses reveal that the learning capacity can also be larger than the equipartition limit at finite sample size. In the Appendix (\ref{sec:appuni}), we work an exactly tractable but non-regular model where the learning capacity is larger than the equipartition value for all sample sizes.

\subsection{The Gibbs entropy and prior selection}

\label{sec:Gibbs}
In statistical physics, the density of states is known (\textit{i.e.}~measured) but in Bayesian inference the selection of a prior is often subjective. The construction of an objective or uninformative prior is a long-standing problem in Bayesian statistics.
What insight does the proposed correspondence provide for prior choice? 
 
Prior construction since Bayes and Laplace has often attempted to apply a \textit{Principle of Indifference}: All \textit{mutually exclusive} and \textit{exhaustive} possibilities should be assigned equal prior probability \cite{Laplace,Keynes1921}. One interpretation of this prescription is that it maximizes entropy  \cite{Jaynes1957information,Shore1980Axiomatic}.
However, the principle of indifference is difficult to interpret in the context of continuous parameters, or across models of different dimension.
For example, are normal models with means $\mu$ and $\mu+{\rm d}\mu$ mutually exclusive (distinguishable)? Even if the mean were constrained to be an integer ($\mu \in \mathbb{Z}$), which would define \textit{mutually exclusive}, the exhaustive condition is also problematic. Exhaustive would correspond to a prior with uniform weighting over all integers. This vanishing prior weight ($1/\infty$) on the non-compact set  $\mathbb{Z}$ results in a paradoxical value for the evidence $Z\rightarrow 0$ and the rejection of the model irrespective of the data,  as described in the Lindley-Bartlett Paradox (Sec.~\ref{sec:lindbartintro}).

\subsubsection{A generalized principle of indifference}
To define \textit{mutually exclusive} in a statistical context, we look for natural analogues to this problem in statistical physics.
%
A surprising result from the perspective of classical physics  is that Nature makes no distinction between states with identical particles exchanged (\textit{e.g.}~electrons) and counts only distinguishable states (the Gibbs paradox).
Following Balasubramanian \cite{Balasubramanian1997}, we proposed that the concept of indistinguishability must be applied to objective Bayesian inference. 
We take the \textit{mutually exclusive} criteria in the principle of indifference to refer to distributions which are mutually distinguishable at the experimental resolution available at sample size $N$.
We propose  a \textit{generalized principle of indifference}: sets of indistinguishable models are each collectively assigned the weight of a single distinguishable model.

To study the weighting of each model, \textit{we must prepare the data using a different procedure}. We distribute $X^N$ according to an assumed true parameter ${\bm \theta}$: $X^N\sim q(\cdot|{\bm \theta})$, omitting the expectation over $\bm \theta$:
\begin{align}
\overline{f}({\bm \theta},N) \equiv \left<\right. f(X^N,{\bm \theta}) \left.\right>_{X},
\end{align}
where $X\sim q(\cdot|{\bm \theta})$. 
A generalized principle of indifference states that the prior $\varpi$ should be chosen such that:
\begin{equation}
\overline{S}({\bm \theta}; N,\varpi) \approx {\rm const}\ \ \ \forall \ \ \ \ {\bm \theta} \in \Theta, \label{Eqn:pofi}
\end{equation}
 at sample size $N$, where the Gibbs entropy is now a function of  ${\bm \theta}$.
Qualitatively, Eq.~\ref{Eqn:pofi} realizes the condition of equal weighting on \textit{mutually exclusive models} since the Gibbs entropy is understood as the log-number of accessible models and constant entropy implies equal weighting between models.



The correspondence \textit{also} offers a natural mechanism for resolving statistical anomalies arising from the \textit{exhaustive} condition in the principle of indifference.
In statistical mechanics, the partition function $Z$ is not a probability  by construction since the density of states $\rho$ is a density but not a probability density. Therefore, a natural solution 
to statistical anomalies arising from the exhaustive condition is to re-interpret the objective inference prior as a \textit{density of models}. 

To circumvent the Lindley-Bartlett Paradox, we must specify a consistent density of models between different parameter values and model families. We replace the  prior $\varpi({\bm \theta})$ with a model density $w({\bm \theta})$ such that: 
\begin{equation}
\overline{S}({\bm \theta};N,w) \approx 0, \label{eqn:gpoi}
\end{equation}
assigning unit multiplicity to all parameters ${\bm \theta}$ and model families $I$. Eq.~\ref{eqn:gpoi} is reparametrization invariant. See the Appendix (\ref{app:invariant}). The prior $w$ will be improper, but none-the-less the normalization is well defined. We shall refer to Eq.~\ref{eqn:gpoi} as the \textit{Generalized Principle of Indifference} which realizes both the mutually-exclusive and exhaustive conditions using a principled statistical approach, regardless of the nature of the parameter manifold. We will call the prior $w$ that satisfies Eq.~\ref{eqn:gpoi} the \textit{GPI prior}.

\subsubsection{Technical note} 
Eqs.~\ref{Eqn:pofi} and \ref{eqn:gpoi} cannot be defined as equalities since the condition is typically not exactly realizable for all ${\bm \theta}$ at finite sample size $N$. To define the GPI prior precisely, we minimize the largest violation of the GPI condition (Eq.~\ref{eqn:gpoi}), using a mini-max approach analogous to that of Kashyap \cite{Kashyap1971}: We  choose the prior $\varpi$ normalization such that 
\begin{equation}
\max_{\bm \theta} \overline{S}({\bm \theta};N,\varpi) = 0, \label{eqn:gpoi2}
\end{equation}
and then the GPI prior is the prior maximizes the minimum Gibbs entropy:
\begin{equation}
w = \underset{\varpi}{\arg\!\,\max} \ \min_{\bm \theta} \overline{S}({\bm \theta};N,\varpi). \label{eqn:gpoi3}
\end{equation}
Qualitatively, this procedure enforces Eq.~$\ref{eqn:gpoi}$ as precisely as possible.

\begin{table*}
\setlength{\tabcolsep}{1em} 
{\renewcommand{\arraystretch}{1.5}
\Scale[0.87]{
 \begin{tabular}{|l|c|c|c|c|c|}
\hline
Model name       &  Likelihood          & Parameters      & Support        & GPI prior & Effective \\
         &  $q(x|{\bm \theta})$ & ${\bm \theta}$  & $\bm \Theta$   & $w({\bm \theta})$  & complexity: $\cal K$ \\
\hline
\hline
Normal  &   $(2\pi \sigma^2)^{-D/2}\ \exp [-\textstyle\frac{1}{2\sigma^2}(\vec{x}-\vec{\mu}\,)^2]$
                  &   $(\vec{\mu})$              &     $\vec{\mu}\in {\mathbb{R}}^D$           & $(\textstyle\frac{N}{2\pi \sigma^2})^{D/2}e^{-{\cal K}}$  & $\textstyle\frac{D}{2}\left[1+N\log(1+N^{-1})\right]$ \\
   &             &  $(\sigma)$ &  $\sigma \in {\mathbb{R}}_+$                   & $(\textstyle\frac{N}{\pi \sigma^2})^{1/2}e^{-{\cal K}}$  & Eq.~\ref{eqn:emum}, $\alpha = \textstyle\frac{1}{2}$ \\                        
             &     & $(\vec{\mu},\sigma)$     &     $\vec{\mu}\in {\mathbb{R}}^D$, $\sigma \in {\mathbb{R}}_+$ & 
             $\sqrt{2}(\textstyle\frac{N}{2\pi \sigma^2})^{(D+1)/2}e^{-\cal K}$   &  Eq.~\ref{eqn:nmumuv} \\ 
             \hline
Exponential  &   $\lambda e^{-\lambda x}$
                  &   $(\lambda) $              &     $\lambda \in {\mathbb{R}}_+$          & $(\textstyle\frac{N }{2\pi\lambda^2})^{1/2}e^{-\cal K}$  &  Eq.~\ref{eqn:emum}, $\alpha = 1$ \\   
                  \hline   
Uniform      &   $H_{\rm SF}(x)H_{\rm SF}(L-x)$
                  &   $(L) $              &     $L \in {\mathbb{R}}_+$          &  $\textstyle\frac{N}{L} e^{-\cal K}$ &    $N\log (1+N^{-1}) + 1$ \\ 
                  \hline
                  Gamma      &   $\textstyle \frac{\beta^\alpha}{\Gamma(\alpha)}x^{\alpha-1}e^{-\beta x}$
                  &   $(\beta) $              &     $\beta \in {\mathbb{R}}_+$          &  $(\textstyle\frac{\alpha N}{2\pi \beta^2})^{1/2} e^{-\cal K}$ &    Eq.~\ref{eqn:emum} \\                    
\hline
\end{tabular} }
\caption{ A summary of exact GPI priors for a collection of standard models. $H_{\rm SF}$ is the Heaviside step function and we write the set of positive real numbers: $\mathbb{R}_{>0}\equiv\{y\in \mathbb{R}|y>0\}$. Many of the effective complexities are defined in the Appendix. \label{tab:GPItab}}}
\end{table*}

\subsection{Examples of the GPI prior}

\label{sec:appGPI}

GPI has properties that rectify significant shortcomings with other approaches to prior selection. Our first aim in this section is to compute the GPI prior for a series of models to show that the calculation is tractable in many applications.
In Sec.~\ref{sec:regmodels}, we compute the GPI prior for regular models in the large-sample-size limit. This analysis reveals a connection between the GPI prior and the Jeffreys prior.
In Sec.~\ref{sec:exact}, we demonstrate an exact computation of the GPI prior for a number of non-harmonic models.


\subsubsection{Approximate GPI prior for regular models}
\label{sec:regmodels}

We will first explore the properties of the generalized principle of indifference by computing the GPI prior in the large-sample-size limit of a regular model.
To define the GPI prior, it is first useful to define the scaled-Jeffreys prior:
\begin{equation}
\rho({\bm \theta};N) \equiv \left(\textstyle\frac{N}{2\pi}\right)^{K/2} I^{1/2}, \label{eqn:scaledJeff}
\end{equation}
where $I$ is the determinant of the Fisher information matrix defined for a single sample (Eq.~\ref{fisherMat}), $K$ is the dimension of the continuous parameter manifold $\Theta$. The prior $\rho$ is a density on the parameter manifold with the qualitative meaning of the inverse volume of indistinguishable models  at sample size $N$. The GPI prior is 
\begin{equation}
\log w({\bm \theta};N) = \log \rho\ -K + {\cal O}(N^{-1}), \label{eqn:w}
\end{equation}
where $K=\dim \Theta$, as shown in the Appendix (\ref{sec:jeffreys}). 

In the large-sample-size limit, the parameter dependence of the GPI prior is identical to the Jeffreys prior, which has enjoyed a long and successful history \cite{Kass1996}. The Jeffreys prior was initially proposed because it was reparametrization invariant \cite{Jeffreys1946}. More recently the same prior has been motivated by numerous other arguments (\textit{e.g.}~\cite{Bernardo1999,Balasubramanian1997}). From the perspective of parameter inference the GPI approach simply recapitulates a widely-applied method in the large-sample-size limit of a regular model. 


\label{sec:nonExact}
\begin{figure}
  \centering
    \includegraphics[width=.5\textwidth]{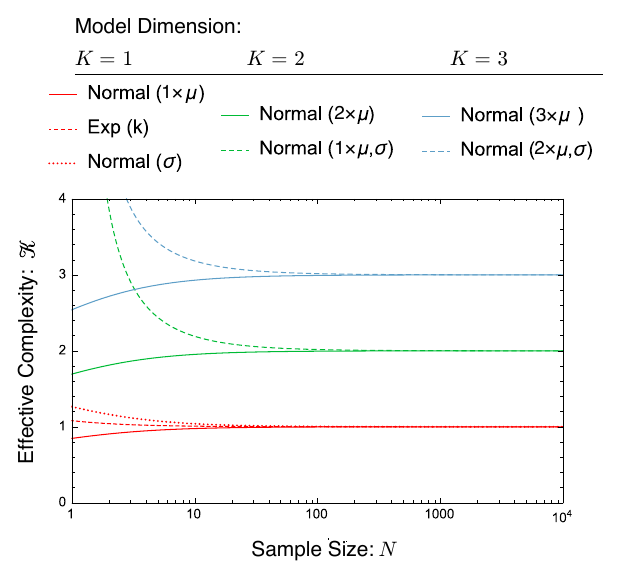} 
    \caption{{\bf Effective complexity of models at finite sample size.} We computed the exact GPI prior for a series of models of different dimension. At large sample size, the dimension determines the effective complexity: ${\cal K} = \dim \Theta / 2$. At finite sample size there are significant corrections. The effective complexity divergences for the normal model (dashed curves) with unknown mean and variance at $N=1$.
      \label{Fig:effectComp}} 
\end{figure}

\subsubsection{Exact GPI prior for models with symmetry}

\label{sec:exact}

For simple models, symmetry and dimensional analysis often imply that $w$ must still be proportional to the Jeffreys prior even at small sample size. We compute the exact GPI prior analytically for the normal model with unknown mean and variance, the exponential model, the uniform model and the Gamma model in the Appendix (\ref{sec:appExactE}). In the Appendix (\ref{sec:appuni}), we demonstrate an exact computation of the GPI prior for a non-regular model. A summary of the exact GPI priors is shown in Tab.~\ref{tab:GPItab}.

All these models have a log-likelihood that is anharmonic in the parameters and therefore are expected to have non-trivial high-temperature behavior.
The calculation reveals that the asymptotic form of the GPI prior (Eq.~\ref{eqn:w}) closely approximates the exact prior. In many models it is convenient to define the finite sample-size correction as an effective complexity $\cal K$ that replaces the model dimension $K$ in Eq.~\ref{eqn:w}:
\begin{equation}
w({\bm \theta};N) = \rho\ e^{-{\cal K}}, \label{eqn:calK}
\end{equation}
$\cal K$ is plotted as a function of sample size ($N$) for a number of different models in Fig.~\ref{Fig:effectComp}. On an empirical basis, it is clear that Eq.~\ref{eqn:w} is typically an excellent approximation for $w$ even small to intermediate sample sizes. 
%


\subsection{GPI circumvents the Lindley-Bartlett paradox}

\label{sec:lind}

%
%

\label{sec:inferenceapp}

To demonstrate that the GPI prior automatically leads to non-anomalous inference (\textit{i.e.}~free from infinite normalization factors) and is also free from \textit{ad hoc} parameters, we return to the example we used to introduce the 
Lindley-Bartlett paradox in the introduction in Sec.~\ref{sec:lindbartintro}: 
We generate data from five competing models, then perform inference on the model identity on each of the datasets using the five models as candidates. A detailed description of the generative parameters is provided in the Appendix (\ref{Sec:CanBayes}).


\subsubsection{GPI approach}

We have computed the GPI prior for each of the proposed models. Inference on parameter values follows the standard Bayesian framework using the GPI prior $w$. 
The GPI prior includes the model prior (Appendix Sec.~\ref{sec:BayesReint}) and therefore the posterior probability of model $I$ is:
\begin{equation}
\varpi(I|x^N) = Z_I/\sum_JZ_J,
\end{equation}
where the model index $J$ runs over the five competing models. The model posteriors for the five sets of simulated data for a sample size of $N = 20$ are shown in Fig.~\ref{Fig:inference}C. The results show a number of important characteristics of the GPI prior: (i) There is an unambiguous  Bayesian procedure for computing inference on both parameters and models. The approach is automatic or free from \textit{ad hoc} or subjective choices. (ii) Inference on both parameters and models leads to non-anomalous results in which the generative distribution has non-zero posterior probability. In our example, the highest posterior model is the generative distribution in each example. (iii) For the normal models, the higher-dimensional models have lower posterior probability for the data generated by model ${\cal N}$, even though the generative distribution is realizable in ${\cal N}(\mu)$ and ${\cal N}(\mu,\sigma)$. This shows that the GPI prior contains an endogenous model selection mechanism favoring model parsimony, and we will discuss this in detail in Sec.~\ref{Sec:ModelSelection}. (iv) Finally, we note that the Jeffreys prior approach is not even possible in the context of the uniform model since the Fisher Information Matrix is undefined. However, we demonstrate in the Appendix (\ref{sec:appExactE}) that the GPI prior can be  computed analytically.


\section{Discussion}

\subsection{Learning capacity}

One valuable feature of the proposed correspondence is the potential to gain new insights into statistical phenomenology using  physical insights into the thermodynamic properties of physical systems. 
Artificial Neural Networks (ANN) and systems-biology models are two examples of systems with a large number of poorly-specified parameters that none-the-less prove qualitatively predictive.  This phenomenon has been termed  \textit{model sloppiness} \cite{Machta:2013hl,Transtrum:2015ud}.  These models often have a logarithmic distribution of Fisher information matrix eigenvalues and this characteristic has been used as a definition of sloppiness  \cite{Machta:2013hl}. But, this definition is unsatisfactory since it is not reparametrization invariant. It is easy to construct counter examples for this definition: For instance, in a $K$-dimensional normal model where the variance for each dimension is logarithmically distributed, the Fisher information eigenvalues are likewise logarithmically distributed, but the model none-the-less behaves like a normal regular model from the standpoint of prediction and statistical analyses.

 The correspondence we describe suggests a definition directly written in terms of the predictive performance of the model and the equipartition theorem. We propose that  \textit{predictive sloppiness} be defined as models that have a smaller learning capacity than estimated from the model dimension:
\begin{equation}
\overline{C} < \textstyle \frac{1}{2} \dim  {\Theta}.
\end{equation}
This definition (i) would exclude all regular models in the large sample-size limit, (ii) is reparametrization invariant and (iii) can be generalized to other non-Bayesian frameworks by expressing the learning capacity in terms of the predictive performance.

\begin{figure}
  \centering
    \includegraphics[width=.48\textwidth]{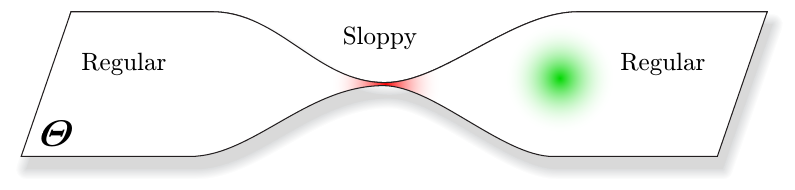}
      \caption{{\bf Sloppiness is determined by parameter manifold geometry and posterior width.} Parameters are defined on a compact manifold ${\bm \Theta}$. In sloppy regions of parameter manifold, the parameters are model-structure dominated (red posterior) whereas in regular regions of parameter manifold parameters are data dominated (green posterior).  From the perspective of the learning capacity, the model is effectively one dimensional in proximity to the red posterior and two dimensional in proximity to the green posterior.
\label{Fig:GPOIA}}
\end{figure}

The mechanism of sloppiness is model parameters being determined by model structure rather than the data. 
This scenario is drawn schematically in Fig.~\ref{Fig:GPOIA}. We sketched a compact parameter manifold after re-parameterizing the model so that the Fisher-Rao metric is the identity matrix. At a regular point (green),  all parameter coordinates are regular and data dominated, the effective dimension of the model is $K_{\rm eff} = 2$. At the sloppy point (red), the manifold is not rigorously singular but the manifold constraints determine the parameter value in the vertical coordinate direction. Therefore the effective dimension is $K_{\rm eff} \approx 1$. In summary, when model structure not data determines the parameter values, the learning capacity will be anomalously small and the model will be  anomalously predictive.

\subsection{The generalized principle of indifference}
\label{sec:GPIdis}

\begin{figure}
  \centering
    \includegraphics[width=.48\textwidth]{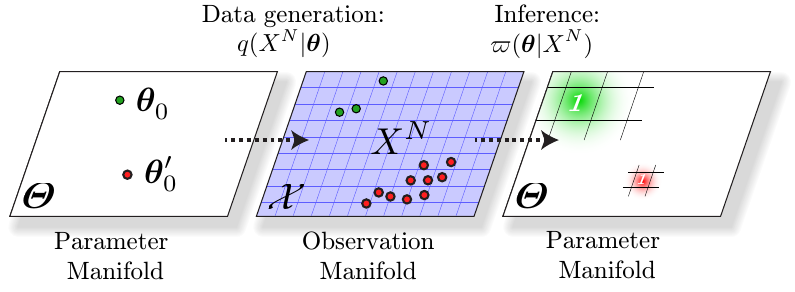}
      \caption{{\bf Generalized principle of indifference.}  The posterior distribution $\inliner{\varpi({\bm \theta}_0|X^N)}$  is shown schematically for two different sample sizes $N$. The resolution increases with sample size  as the posterior shrinks.  In the GPI  prior (Eq.~\ref{eqn:gpoi}), all parameter values consistent with the posterior are assigned  unit prior weight collectively. 
\label{Fig:GPOIB}}
\end{figure}

We argue that a natural approach to objective Bayesian inference is to choose a  prior  such that the number of indistinguishable distributions is one for all parameter values. (See Eq.~\ref{eqn:gpoi}.) Schematically, this procedure assigns equal prior weighting to all models that can be distinguished at finite sample size $N$. As the sample size increases, the prior must be modified to accommodate the increased resolution (Eq.~\ref{eqn:resolution}) due to shrinking of the posterior support.  (See Fig.~\ref{Fig:GPOIB}.) For a regular model in the large-sample-size limit, no calculation is required and GPI  prior is equal to the scaled-Jeffreys prior (Eq.~\ref{eqn:scaledJeff}.)  At small sample size or in singular models, the GPI prior must be computed explicitly.  

It is important to stress that GPI gives rise to a \textit{sample-size-dependent prior} and therefore this inference is \textit{not} Bayesian in a classical sense: (i) It violates Lindley's dictum: \textit{today's posterior is  tomorrow's prior}. (ii) Furthermore, the evidence and prior are no longer interpretable as probabilities but rather statistical weightings. On-the-other-hand, the method codes parameter uncertainty in terms of a posterior probability distribution and facilitates Bayesian parameter and model averaging. Therefore, we would argue the approach maintains  all of the  attractive features of the Bayesian framework while avoiding problematic aspects.   

\subsubsection{Model selection}
\label{Sec:ModelSelection}
The normalization of the GPI prior has  significant consequences for inference on model identity (\textit{i.e.}~model selection). Returning to the regular model, it is straight forward to apply the Laplace approximation to compute the minus-log evidence using the GPI prior:
\begin{equation}
-\log Z(x^N;w) \approx -\log q(x^N|\hat{\bm \theta}) + K, \label{eqn:mle}
\end{equation}
where $K=\dim \Theta$. The scaled-Jeffreys prior cancels the Occam factor from the integration. The two remaining contributions each have clear qualitative interpretations: the MLE estimate of the information ($-\log q$) is a measure of the \textit{goodness-of-fit} and a penalty for model complexity  ($K$).  Eq.~\ref{eqn:mle} is already well known as the Akaike Information Criterion (AIC)\footnote{Where AIC is defined in nats, rather than the more common demi-nat expression which is twice Eq.~\ref{eqn:mle}.}:
\begin{equation}
{\rm AIC}(x^N) \equiv -\log q(x^N|\hat{\bm \theta}) + K. \label{eqn:aic}
\end{equation}
Information-based inference is performed by selecting the model which minimizes AIC, maximizing the estimated predictive performance \cite{BurnhamBook}. The reason why AIC and GPI Bayesian inference are equivalent is most easily understood by rewriting Eq.~\ref{eqn:gpoi}:
\begin{equation}
-N\overline{F} \approx -N\overline{U}, \label{eqn:LSSL}
\end{equation}
which in statistical language corresponds to using a prior that makes the log partition function (LHS) an unbiased estimator of the log predictive performance (RHS). 
Since the Akaike Information Criterion (AIC) is an unbiased estimator of RHS at large sample size $N$, the generalized principle of indifference encodes an AIC-like model selection \cite{Stone1977} and an information-based (AIC) realization of Occam's razor: \textit{parsimony increases predictivity} \cite{BurnhamBook}.
The log-predictive performance (RHS Eq.~\ref{eqn:LSSL}) has already been advocated in the context of Bayesian model selection through the use of pseudo-Bayes factors by Gelman and coworkers \cite{Gelman2014,gelfand1994bayesian,Spiegelhalter2002,Vehtari2012,BURNHAM2004}.

\subsubsection{Posterior impropriety}

\label{sec:improp}

The use of GPI prior often, but not always, gives non-zero evidence for all models under consideration.
One such exception is shown in Fig.~\ref{Fig:effectComp} which reveals that the normal model with unknown mean and variance has a divergent effective complexity at a sample size of $N=1$. The effect of this divergence is to give these models zero statistical weight. Although this may initially appear problematic, it is an important feature of the generalized principle of indifference. A mean and variance cannot be estimated from a single observation and as a result the model parameter posterior would be improper and the predictive loss would be infinite. The generalized principle of indifference automatically removes these models from consideration by assigning these models  zero statistical weight. The inability  to automatically handle posterior impropriety is recognized as a significant shortcoming  of these competing approaches \cite{Kass1996}.

\subsection{Comparison with existing approaches}
The generalized principle of indifference subsumes a patchwork of conflicting methods for prior and model selection, resolving many conflicting approaches and generating a single,  generally-applicable and self-consistent framework.
The GPI approach subsumes the following approaches: (i) For discrete parameter manifolds in the large sample size limit, the GPI gives equal weight to all mutually exclusive models, consistent with the original formulation of the principle of indifference by Bayes and Laplace  \cite{Laplace,Keynes1921}. (See Eq.~\ref{eqn:discretew}.) (ii) In the large-sample-size limit, GPI generates a GPI prior proportional to the well-known Jeffreys prior. In this sense, the approach is closely related to the reference prior approach of Bernardo and Berger \cite{BernardoBerger1991,Bernardo1999}. (iii) With respect to model selection (inference on model identity), the GPI evidence behaves like pseudo-Bayes factors (or AIC) and therefore circumvents the Lindley-Bartlett paradox. (See Sec.~\ref{Sec:ModelSelection}.) To date, the pseudo-Bayes approach has always been un-Bayesian in the sense that 
the pseudo-Bayes method consists of the \textit{ad hoc} combination of a canonical Bayesian prior for inference on parameters but a cross-validation-based weighting for inference on models. The GPI provides a self-consistent  approach to inference on both parameters and models.

The GPI addresses a number of problems with existing approaches to objective Bayesian inference.
(iv) \textit{Lindley-Bartlett paradox:} As already discussed above,  an important shortcoming with existing objective Bayesian approaches relates to the compactness of the parameter manifold and the automatic rejection of higher-dimensional models in model selection (the Bartlett-Lindley paradox \cite{Lindley1957,Bartlett1957}).  More generally, the evidence of the canonical objective Bayesian approach depends on \textit{ad hoc} modeling decisions, like the range of allowed parameter values.
The GPI-Bayes evidence circumvents these anomalies by generating a consistent distribution  density $w$ over competing models. As a result the GPI evidence is independent of \textit{ad hoc} modeling decisions.  (v) \textit{Unification of statistical paradigms:} The absence of the Lindley-Bartlett paradox implies coherent inference between paradigms \cite{Cousins2014,LaMont2017b} and therefore the generalized principle of indifference naturally unifies objective Bayesian inference  with information-based inference.
(vi) \textit{Prior and posterior impropriety:} Another important flaw identified in other objective Bayesian approaches is the inability to handle impropriety. In many cases where parameters are defined on non-compact manifolds, the prior (and sometimes the posterior) cannot be normalized. The redefinition of the prior as a density of models introduces a well-defined and consistent method for defining prior normalization, regardless of the global structure of the manifold. Furthermore, the approach automatically assigns zero statistical weight to models that suffer from posterior impropriety. (See Sec.~\ref{sec:improp}.) 
(vii)
\textit{Singularity and sloppiness:} Finally, the GPI-Bayes approach does not assume model regularity. It treats singularity and the sloppiness phenomenon in a natural way. (See Sec.~\ref{sec:exact}.)

\subsection{Conclusion}

Nature reveals an elegant formulation of statistics in the thermal properties of physical systems. Measurements of the heat capacity, compressibility or susceptibility reveal unambiguously how Nature enumerates states and defines entropy.
These physical insights provide clues to the definition of novel statistical quantities and the resolutions of ambiguities in the formulation of objective Bayesian statistics.
We have refined a previously proposed correspondence between the Bayesian marginal likelihood and the  partition function of statistical physics. 
We demonstrate a novel and substantive mapping between the average energy, heat capacity, entropy and other statistical quantities. 
The newly-defined learning capacity is a natural quantity for characterizing and understanding  learning algorithms and generates new physical insight into  the mechanism of model sloppiness through a correspondence with the equipartition theorem and the freeze-out phenomenon. 
A key motivation for exploring the phenomenology of learning is to apply these insights to develop new learning algorithms. We provide one example in the paper: We use the Gibbs entropy to define a generalized principle of indifference and an objective Bayesian GPI prior with the property that all distributions have equal prior weight.
This approach subsumes many seemingly inconsistent and disparate methods into a single, coherent statistical approach.  For the first time, we demonstrate a self-consistent Bayesian approach to performing inference on models of unknown dimensional with uninformative priors.

\medskip
\noindent
\textbf{Acknowledgements:} PAW and CHL acknowledge helpful discussions with M.~Linden, J.~Kinney, D.~Mayo, C.~Heilig, M.~Abbott, B.~Machta and M.~Transtrum. This work was support by NSF grants NSF-PHY-084845 and NSF-MCB-1151043-CAREER.

\newpage
\onecolumngrid
\appendix


\section{Supplemental results}

\subsection{Details on the Lindley-Bartlett Paradox example}

\label{Sec:CanBayes}

\subsubsection{Canonical objective Bayesian approach}

\noindent
\textit{Analysis:}
For the canonical objective Bayesian approach, we will use a proper objective prior. We attempt to use the Jeffreys prior for each model.
A problem immediately presents itself in the context of the Uniform model where the Jeffreys prior is undefined.  (See Sec.~\ref{sec:appuni}.)
We must therefore deviate from our protocol and apply some other prior. We set a flat prior on the parameter $L$ motivated by the principle of indifference.  The priors for the four models with parameters cannot be normalized due to their non-compact parameter manifolds. Formally, we can work in a finite interval, then consider the limit as the limits of the intervals approach infinity.

\smallskip
\noindent
\textit{Parameters:} Parameter posteriors for each model can be computed using this procedure and the results are identical to the GPI approach, except for the uniform model where no Jeffreys prior exists. These are clearly acceptable results.

\smallskip
\noindent
\textit{Models:}  Unlike the parameter posteriors,  inference on the model leads to anomalous results. We find that the model posterior for the parameter-free normal model is one, regardless of which distribution was used to generate the data, due to the prior impropriety of the other models. (See Fig.~\ref{Fig:inference}A.)

\subsubsection{Revised uninformative Bayesian approach}

After having seen the data, a Bayesian will often reconsider the prior and localize it around the values favored by the data. Here we normalize the priors on the revised finite intervals defined in Tab.~\ref{tab:models}. It is important to stress that the boundary of each interval is \textit{ad hoc} and investigators will make different choices. Sensible choices typically do not strongly affect the parameter posteriors, but they do affect model posteriors as shown in Fig.~\ref{Fig:inference}B. This approach is formalized in variational or empirical Bayesian methods. In this case, the prior is no longer determined \textit{a priori} and this double-use of data can lead to difficulties due to the potential for overfitting.

\begin{table}
\setlength{\tabcolsep}{1em} 
{\renewcommand{\arraystretch}{1.5}
\Scale[0.87]{
  \begin{tabular}{|lc|c|c|c|c|c|}
\hline
Model & & Initial parameter  & Revised parameter  & Generative parameters & Likelihood & Prior \\
      & & support: ${\bm \theta}$   & support: ${\bm \theta}$  &   ${\bm \theta}_0$ & $q(x|{\bm \theta} )$ & $\varpi({\bm \theta})$ \\ 
\hline
\hline
Normal & ${\cal N}$ &  & & $\mu_0 = 5$, $\sigma_0 = 1$   & $(2 \pi \sigma^2)^{-1/2}e^{-(x-\mu)^2/2\sigma^2}$ & 1 \\
 & ${\cal N}(\mu)$ & $\mu\in {\mathbb R}$ &  $\mu\in [0,10]$ & $\mu_0 = 6$, $\sigma_0 = 1$  &  &  $C_0$ \\
& ${\cal N}(\mu,\sigma)$ & $\mu\in {\mathbb R}$, $\sigma\in {\mathbb R}_{>0}$ & $\mu\in [0,10]$, $\sigma\in [0.1,10]$ & $\mu_0 = 5$, $\sigma_0 = 0.75$  &  & $C_0/\sigma^2$  \\
\hline
Exponential & ${\rm Exp}(\lambda)$ & $\lambda\in {\mathbb R}_{>0}$ & $\lambda\in [0.1,10]$ & $\lambda_0 = 2$  & $\lambda e^{-\lambda x}$ & $C_0/ \lambda$ \\
\hline
Uniform & ${\cal U}(L)$ & $L\in {\mathbb R}_{>0}$ & $L\in [0,10]$ & $L_0=10$ & $H_{\rm SF}(x)H_{\rm SF}(L-x)/L$ & $C_0$ \\
\hline
\end{tabular} }
\caption{{\bf Models for inference on simulated data.}  Five datasets were generated, one for each model, using the generative parameters: $\inliner{X^N\sim q(\cdot|{\bm \theta}_0)}$. In the canonical objective Bayesian approach, inference was performed on the simulated data using the the objective prior shown in the final column. The normalization constant $C_0$ was chosen in each case to make the prior proper over the original or revised parameter support. In the GPI approach, the relevant GPI prior $w$ was used, as computed in this appendix, for each model. See Tab.~\ref{tab:GPItab}. $H_{\rm SF}(x)$ is the Heaviside step function. \label{tab:models}}}
\end{table}

\subsection{Definitions of information, cross entropy, Fisher information matrix}

\label{sec:intro}

The Shannon Information is defined:
\begin{equation}
h(x|{\bm \theta}) \equiv -\log q(x|{\bm \theta}).
\end{equation}
Let $X$ be distributed with a true distribution with parameter ${\bm \theta_0}$: $X \sim q(\cdot|{\bm \theta_0})$. The cross entropy is defined:
\begin{equation}
H({\bm \theta};{\bm \theta}_0) \equiv \left< h(X|{\bm \theta})\right>_{X\sim q(\cdot|{\bm \theta_0})},  \label{eqn:Ce}
\end{equation} 
and which has a minimum at the true entropy:  
\begin{equation}
H_0({\bm \theta}_0) \equiv H({\bm \theta}_0;{\bm \theta}_0).
\end{equation} 
The empirical estimator of the cross entropy is defined:
\begin{equation}
\hat{H}({\bm \theta}) \equiv N^{-1} \sum^{N}_{i=1} h(x_i|{\bm \theta}),\label{eqn:emp}
\end{equation}
which is independent of $N$ to leading order, in spite of the prefactor.
The KL-Divergence:
\begin{equation}
D_{\rm KL}({\bm \theta}_0||{\bm \theta}) = H({\bm \theta};{\bm \theta}_0)-H_0({\bm \theta}_0),
\end{equation}
is the natural distance-like measure on the parameter manifold. The Fisher information matrix is defined:
\begin{equation}
\label{fisherMat}
I_{ij} = \left[ \textstyle \frac{\partial}{\partial \theta^i}\textstyle \frac{\partial}{\partial \theta^j} H({\bm \theta};{\bm \theta}_0) \right]_{\bm \theta=\bm \theta_0}, 
\end{equation}
which is a rank-two symmetric covariant tensor known as the Fisher-Rao metric \cite{Amari1985}.

\subsection{An alternate correspondence}
\label{App:alttemp}

We do not believe that statistical mechanics prescribes a unique procedure for objective Bayesian inference.
In establishing the correspondence between inference and statistical mechanics, we identify the partition function $Z$ as the marginal likelihood and $N \leftrightarrow \beta$ in agreement with V.~Balasubramanian \cite{Balasubramanian1997}. However, this is not the only choice that has been proposed. For instance, Watanabe \cite{watanabe2009} instead chooses to define  the inverse-temperature $\beta^*$ so that the likelihood is raised to an arbitrary power $\beta^*$:
\begin{equation}
q(X^N|{\bm \theta}) \rightarrow q^{\beta^*}(X^N|{\bm \theta}). 
\end{equation}
We have explored both possibilities in some detail. In the large-sample-size limit, these two correspondences have many similar properties. However, we present only the analysis of the $N \leftrightarrow \beta$ correspondence in this papers because we feel that it identifies statistical objects with the most desirable properties.

It is important to note that the $\beta^*$ correspondence does have a number of convenient computational properties:
 (i) It is continuous and therefore no finite difference definition need be introduced. (ii) It allows one to interpolate between a Bayesian posterior (given by $\beta^* = 1$) and the point estimates of the MLE's (given by $\beta^* \rightarrow \infty$). This temperature has also been applied in tempering schemes in MCMC methods, and simulated annealing---increasing the temperature promotes a better exploration of the sample space (chain-mixing) that can be used to better sample multimodal distributions, or find the minima in a rough function in close analogy to analogous methods in physics.

However,  the $\beta^*$ correspondence also has a number of features we find undesirable: First, $\beta^*$ is not a preexisting statistical parameter within the Bayesian framework. Only $\beta^*=1$ corresponds to a Bayesian statistics. High and low $\beta^*$ correspond to non-Bayesian statistics whereas the sample size interpretation of $\beta$ has a sensible and natural Bayesian interpretation in terms of small and large sample sizes. 
Second, the internal energy under the choice of $\beta^*$ is not the predictive performance $U$. Consequently, the principle of indifference, which results from a likelihood-power $\beta^*$, does not generate the Akaike weights as the model averaging procedure. 
 It is the reproduction of  AIC, with its proven asymptotic efficiency \cite{Shao1997}, that is an important motivation of the proposed correspondence.

\subsection{Finite difference is equivalent to cross validation}
\label{App:crossval}
The log-predictive distribution can be written as a finite difference:
\begin{align}
\log q(X_i|X^{\ne i}) &= \log Z(X^{N}) - \log Z(X^{\ne i}).
\end{align}
We can interpret the $-\log q(X_i|X^{\ne i})$ as a finite difference estimate of the the sample size derivative of the free energy. We take the mean over all permutations of the data so that this estimate is symmetric with respect to all data points. Under expectation, analytically continuing sample size, the LOOCV relationship to the internal energy is clear:
\begin{align}
\langle \log q(X_i|X^{\ne i}) \rangle &\approx \frac{\partial}{\partial N} \langle \log Z(X^{N}) \rangle + O(N^{-1}).
\end{align}
This identity is crucial in establishing the thermodynamic interpretation in terms of predictive performance.

\subsection{Jeffreys prior is proportional to GPI prior in the large-sample-size limit}

\label{sec:jeffreys}
In the large-sample-size limit, the partition function can be evaluated using the Lapalce (saddle-point) approximation and the resulting prior is proportional to the Jeffreys prior.  The integral is evaluated by expanding around the minimum of $\hat{H}_X({\bm \theta})$, the maximum likelihood estimator: $\hat{\bm \theta}_X$. The partitition function $Z(X^N) =  \int_{\bm \Theta} \!\! {\rm d}{\bm \theta}\ \varpi({\bm \theta})\,\exp[-N\hat{H}_X({\bm \theta})]$, becomes
\begin{align}
Z(X^N) 
&\approx e^{-N \hat{H}_X(\hat{\bm \theta}_X)} \left(\frac{2 \pi}{N(\det I)^{1/K}}\right)^{K/2}  \varpi({\bm \theta}_X)
\intertext{By the standard $\chi^2_K$ representation of the overfitting error, $\langle \hat{H}(\hat{\bf \theta}_X) \rangle_X = H_0 - \frac{K}{2N}$. Therefore the disorder average becomes}
\overline{F}({\bm \theta}_0,\varpi,N) &=    H_0 - \frac{K}{2N} - \frac{K}{2N} \log \frac{2 \pi}{N(\det I)^{1/K}} -\frac{1}{N} \log \varpi({\bm \theta}_0)+ O(N^{-2})\\
\intertext{We can then calculate the Gibbs entropy $N^2 \partial_N F$,}
\overline{S}({\bm \theta}_0,\varpi,N) &=  \frac{K}{2} \log \frac{2 \pi}{N(\det I)^{1/K}} + K +  \log \varpi({\bm \theta}_0) + O(N^{-1})
\end{align}
If we enforce the generalized principle of indifference, ignoring order $N^{-1}$,
\begin{align}
0 &= S({\bm \theta}_0,w,N)\\
\intertext{and substituting the $w$ for $\varpi$ in the entropy expression then gives us the condition}
w({\bm \theta}_0) &= (\det I)^{1/2} \left( \frac{N}{2\pi} \right) ^{K/2} e^{-K}.
\end{align}
Thus the GPI condition is satisfied by the Jeffries prior in the large-sample-size limit.  The constant weighting factor  is important in model selection as $e^{-K}$ encodes the Akaike weight \cite{BurnhamBook}. The GPI prior has sample-size dependence. 
This sample size dependence will break the de-Finetti likelihood principle: that the prior should not depend on the nature of the data-generating procedure (including the sample size) \cite{Bernardo1994}. The departure from the likelihood principle is the origin of the departure from the conventional Bayesian model selection behavior. 

\subsection{Reparametrization invariance of the GPI approach}
\label{app:invariant}
An important property of an objective prior is \textit{reparameterization invariance}. Is the GPI approach reparameterization invariant? First consider the properties of the partition function:
\begin{align}
Z(X^N) &= \int \mathrm{d}{ \bm \theta} \, \varpi( {\bm \theta} )\, q(X^N|{\bm \theta}), 
\end{align} 
 which is  invariant under reparameterization ${\bm \theta} \rightarrow {\bm \theta}'$ if the prior $\varpi$  transforms as a density, \textit{i.e.}:
\begin{align}
\varpi'({\bm \theta}') = J^{-1} \varpi( {\bm \theta} ),
\end{align}
where $J$ is the determinant of the Jacobian of the coordinate transformation. Since the GPI condition is written in terms of the partition function and the sample size $N$, the GPI condition itself is also reparameterization invariant. Therefore if $w({\bm \theta})$ is the GPI prior for $\bm \theta$ coordinates, 
\begin{align}
w'({\bm \theta}') = J^{-1} w( {\bm \theta} ),
\end{align}
will be the GPI prior in the  $\bm \theta'$ coordinates. Therefore the GPI prior will transform like a density under reparameterization.

\subsection{Effective temperature of confinement}
\label{App:freeparticle}
To calculate the free energy $F$ of a classical free particle confined to a volume $V=L^3$, we calculate the partition function by integrating over available phase space:
\begin{align}
{Z}(\beta) &= \int \frac{d^K {\bm p} \, d^K {\bm x} }{(2 \pi \hbar)^K} e^{-\beta H({\bm p},{\bm x}) } \\
&= \frac{e^{-\beta E_0} L^K}{(2\pi \hbar )^K} \left( \int dp \, e^\frac{-\beta p^2 }{2m}\right )^K = \left( \frac{m L^2}{ 2\pi \hbar^2 \beta} \right)^{K/2}e^{-\beta E_0} .
\end{align}
The Free energy is then
\begin{align}
F(\beta) = E_0 + \frac{K}{2\beta} \log \frac{m L^2}{ 2\pi \hbar^2 \beta} = E_0 + \frac{K}{2\beta} \log \frac{\beta_0}{\beta}
\end{align}
where we have made the identifications
\begin{align}
\beta_0 =  \frac{m L^2}{ 2\pi \hbar^2} \quad \mathrm{and} \quad K = 3.
\end{align}
$\beta_0$ can be interpreted as the inverse of the (typically negligibly small) temperature at which the thermal de Broglie wavelength of the confined particle is on the order of the width of the confining box.

\subsection{ A Bayesian re-interpretation}

\label{sec:BayesReint}

The replacement of the prior (a probability density) with an unnormalized density of states may make a Bayesian reader uncomfortable since the evidence ($Z$) no longer has the meaning of a probability. But there is a natural Bayesian interpretation in terms of the  model prior. 

Typically, when models are compared in a Bayesian context, all  mutually exclusive models are  assigned equal \textit{a priori} probabilities (\textit{i.e.} the principle of indifference). But, we have now proposed a new concept of model enumeration by introducing a density of models. We can compute the total number of distinguishable distributions in model $I$ at sample size $N$ by integrating the GPI prior (density of states) over the parameter manifold: 
\begin{equation}
{\cal N}_I(N) \equiv \int_{\Theta}\!\! {\rm d} {\bm \theta}\ w_I({\bm \theta};N).
\end{equation}
Since models $I$ and $J$  contain different numbers of distinguishable distributions, we reason that the principle of indifference should be interpreted to apply at the distinguishable distribution level rather than the model level. Therefore the \textit{a prior} model probabilities should be:
\begin{equation}
\varpi_I\equiv{\cal N}_I/\textstyle\sum_J{\cal N}_J. \label{eqn:modelNumber}
\end{equation}
and the proper parameter prior is 
\begin{equation}
\varpi({\bm \theta}|I) \equiv w_I({\bm \theta};N)/{\cal N}_I. \label{eqn:parameterPrior}
\end{equation}
Inference with the improper GPI prior is equivalent to assuming proper prior $\varpi_I$ on models and proper prior $\varpi({\bm \theta}|I)$ on parameters.  The numerator in  RHS of Eq.~\ref{eqn:modelNumber} will cancel the denominator in the RHS of Eq.~\ref{eqn:parameterPrior} when the model posterior is computed and the normalization ${\cal N}_I$ divides out of parameter posterior distributions.

%
%

\section{Methods}

\subsection{Computation of the free energy using a sufficient statistic}
\label{SufficiencyCorrection}

It is often convenient to work in terms of sufficient statistics because (i) all the data dependence of the posterior enters through the sufficient statistic and (ii) the statistics have well known statistical distributions that significantly simplify many calculations. We define a sufficient statistic ${\bm t} = {\bm T}(X^N)$ such that 
\begin{equation}
{\rm Pr}({\bm \theta}|X^N) = {\rm Pr}({\bm \theta}|{\bm t}),
\end{equation}
or all the information about the parameters is encoded in ${\bm t}$. We can therefore write:\begin{equation}
q(X^N|{\bm \theta}) = q(X^N|{\bm t})\ q({\bm t}|{\bm \theta}),
\end{equation}
and we can define a Statistic Shanon entropy:
\begin{equation}
H_{\bm t}({\bm \theta}) = -\overline{\log q({\bm t}|{\bm \theta})}.
\end{equation}
In terms of the sufficient statistic, the partition function factors:
\begin{eqnarray}
Z(X^N;\varpi) &=& q(X^N|{\bm t})\ z({\bm t};N,\varpi),
\end{eqnarray}
where the statistic partition function is
\begin{eqnarray}
z({\bm t};N,\varpi) &\equiv& \int_{\Theta}{\rm d}{\bm \theta}\ \varpi({\bm \theta})\ q({\bm t}|{\bm \theta}).
\end{eqnarray}
The expected free energy can  be written:
\begin{equation}
\label{eq:statpartition}
\overline{F}({\bm \theta};N,\varpi) = -N^{-1}\overline{\log z({\bm t};N,\varpi)}+H_0({\bm \theta})-N^{-1}H_{\bm t}({\bm \theta}),
\end{equation}
where $H_0$ is the entropy.



\subsection{Computation of learning capacity}

To compute the learning capacity, we will use the definition from Tab.~\ref{Table:stat}:
\begin{equation}
\overline{C}({\bm \theta};N,\varpi) = N^2\partial_N^2 \overline{\log Z}({\bm \theta};N,\varpi), 
\end{equation}
where $X\sim q(\cdot|{\bm \theta})$. We will ignore the finite-difference definition in these computations for simplicity.

\subsection{Direct computation of GPI prior}
\label{sec:appExact}

We will use the finite-difference definition of the entropy (Eqs.~\ref{eqn:loocv} and \ref{eqn:loocv2}) to enforce the generalized principle of indifference (Eq.~\ref{eqn:gpoi}). The relation for the GPI prior can be written: 
\begin{equation}
(N+1)\,\overline{\log Z}({\bm \theta};N,w) = N\,\overline{\log Z}({\bm \theta};N+1,w), \label{eqn:part} 
\end{equation}
in terms of the partition function. We will use Eq.~\ref{eqn:part} explicitly to solve for the GPI prior $w$.
For the models we work analytically, we will be be able to use the asymptotic form of $w$ (Eq.~\ref{eqn:w}) to define an effective model dimension ${\cal K}$ (Eq.~\ref{eqn:calK}).  The general strategy will be:
\begin{enumerate}
\item{Use symmetry and dimensional analysis to deduce the scaling of $w$ with respect to the parameters ${\bm \theta}$.}
\item{Compute $\log Z(X^N;w)$ and re-express in terms of canonical random variables.}
\item{Compute $\overline{\log Z(X^N;w)}$.}
\item{Solve for the unknown normalization $c$ of $w$ using the GPI condition (Eq.~\ref{eqn:part}).}
\end{enumerate}

\subsection{Computation of the GPI prior using a recursive approximation}
\label{sec:recursiveAlgorithm}
The Gibbs entropy has the property that it is linear in the prior such that:
\begin{align}
\overline{S}(\theta_0, N, e^{\alpha} \varpi) = \alpha + \overline{S}(\theta_0, N, \varpi)
\end{align}
If the prior and entropy are flat, then setting $\alpha = -\overline{S}(\theta_0, N, \varpi)$ will result in $\overline{S}(\theta_0, N, e^{\alpha} \varpi) = 0 $; the  GPI prior condition. This suggests the following simple recursive scheme for a successive approximation for the GPI prior:
\begin{algorithmic}[1]
\Procedure{RecursiveW}{$\varpi$}
\Repeat 
  \State $\varpi(\theta)  \gets \varpi(\theta) e^{-\overline{S}(\theta,\varpi,N)}$
\Until{$\overline{S}(\theta; \varpi, N) \approx 0$}
\EndProcedure
\end{algorithmic}
To the extent the entropy is slowly varying and only locally dependent on the prior, this algorithm will very quickly converge to an exact GPI prior. However, effects due to manifold boundaries and model singularities may create artifacts that lead to unstable updates. Empirical evidence suggests that the algorithm should be terminated before the exact GPI prior condition is met.  Typically very few iterations are required. 
%
%
%
%
%
%

\section{Detailed analysis of the learning capacity and GPI prior of selected models}

\label{sec:appExactE}

\subsection{Normal models}

\label{app:norm}

\subsubsection{Normal model with unknown mean, known variance and an informative prior}

\label{sec:mmumkvp}

The  likelihood for the normal model is defined by Eq.~\ref{eqn:normalmodel}
with parameters ${\bm \theta} \equiv (\vec{\mu})$ for support $\mu \in \mathbb{R}^K$ for a normal model with unknown mean and known variance $\sigma^2$. In this example, we assume a conjugate prior:
\begin{equation}
\varpi({\bm \theta}) = (2 \pi \sigma^2_\varpi)^{-K/2} \exp[ -\textstyle\frac{1}{2\sigma^2_\varpi}(\vec{\mu}-\vec{\mu}_\varpi)^2],
\end{equation}
where we define the critical sample size:
\begin{equation}
N_0 \equiv \sigma^2/\sigma_\varpi^2.
\end{equation}
The partition function is computed by completing the square in the exponential.
If $X^N \sim q(\cdot|{\bm \theta})$ and ${\bm \theta} \sim \varpi$, the log partition function can be expressed in terms of three  chi-squared random variables:
\begin{eqnarray}
\sigma^{-2}\sum_{i=1}^N(\vec{X}_i-\hat{\vec{\mu}}_X)^2 &\sim& \chi^2_{K(N-1)}, \\
\sigma^{-2}N(\vec{\mu}-\hat{\vec{\mu}}_X)^2 &\sim& \chi^2_{K},\\
\sigma^{-2}N_0(\vec{\mu}-\vec{\mu}_\varpi)^2 &\sim& \chi^2_{K}.
\end{eqnarray}
The log partition function is therefore distributed:
\begin{equation}
\log Z( X^N;\varpi) \sim -\textstyle \frac{KN}{2}\log 2 \pi \sigma^2-\textstyle \frac{K}{2}\log\textstyle\frac{N+N_0}{N_0} -\textstyle \frac{1}{2}\chi^2_{K(N-1)} -\textstyle\frac{1}{2}\textstyle\frac{N_0N}{N+N_0}( N^{-1}\chi^2_K+N_0^{-1}{\chi^2_K}'),
\end{equation}
where $\chi^2_j$ is a chi-squared random variable dimension $j$ and the expect-log partition function is 
\begin{equation}
\overline{\log Z}(N,\varpi) = -NH_0-\textstyle \frac{K}{2}\log\textstyle\frac{N+N_0}{N_0},
\end{equation}
where $H_0$ is the entropy and the free energy is:
\begin{equation}
\overline{F}(N,\varpi) = H_0+\textstyle \frac{K}{2N}\log\textstyle\frac{N+N_0}{N_0}.
\end{equation}
The other results in the Tab.~\ref{tab:results} are generated by apply the definitions of the correspondence in Tab.~\ref{Table:stat}.

\subsubsection{Normal model with unknown mean and known variance}

\label{sec:mmumkv}

The  likelihood for the normal model is defined by Eq.~\ref{eqn:normalmodel}
with parameters ${\bm \theta} \equiv (\vec{\mu})$ for support $\vec{\mu}\in \mathbb{R}^D$ for a normal model with unknown mean and known variance $\sigma^2$. The cross entropy is
\begin{equation}
H({\bm \theta};{\bm \theta}_0) \equiv \textstyle\frac{D}{2} \left[\log 2\pi \sigma^2 + 1\right]+\textstyle\frac{1}{2\sigma^2} (\vec{\mu}-\vec{\mu}_0)^2, 
\end{equation}
where the true distribution is $X\sim q(\vec{x}|{\bm \theta}_0)$ and the determinant of the Fisher information matrix is:
\begin{equation}
\det {\bm I} = \sigma^{-2D}.
\end{equation}
The scaled Jeffreys prior (Eq.~\ref{eqn:scaledJeff}) is therefore: 
\begin{equation}
\rho = (\textstyle\frac{N}{2\pi \sigma^2})^{D/2}.
\end{equation}
We will assume $w$ matches the asymptotic form:
\begin{equation}
w = c \sigma^{-D},
\end{equation}
and solve for the unknown constant $c(N,D)$. The partition function is 
\begin{equation}
\log Z(X^N;w) \sim \log c - \textstyle\frac{DN}{2} \log 2\pi \sigma^2 + \textstyle\frac{D}{2}\log \textstyle\frac{2\pi}{N}-\textstyle\frac{1}{2}\chi^2_{D(N-1)}
\end{equation}
where $\chi^2_{D(N-1)}$ is a $D(N-1)$-dimensional chi-squared random variable. The expected log partition function is:
\begin{equation}
\overline{\log Z}({\bm \theta};N,w) = -NH_0({\bm \theta})+ \log c + \textstyle\frac{D}{2}\left[\log \textstyle\frac{2\pi}{N}+1\right],
\end{equation}
where $H_0$ is the true entropy. The learning capacity is:
\begin{equation}
\overline{C}({\bm \theta};N) = \textstyle\frac{D}{2},
\end{equation}
where $D$ is both the dimension of the model. The
unknown normalization of $w$ is:  
\begin{equation}
\log c = \textstyle\frac{D}{2}\log \textstyle\frac{N}{2\pi} - \textstyle\frac{D}{2}\left[1+N\log(1+N^{-1})\right]
\end{equation}
which can be re-written as an effective complexity:
\begin{equation}
{\cal K} =  \textstyle\frac{D}{2}\left[1+N\log(1+N^{-1})\right],
\end{equation}
to define the GPI prior $w$ using Eq.~\ref{eqn:calK}.

\subsubsection{ Normal model with unknown discrete mean}
\label{App:discra}

The  likelihood for the normal model is defined by Eq.~\ref{eqn:normalmodel}
with parameters ${\bm \theta} \equiv (\vec{\mu})$ for support $\vec{\mu}\in \mathbb{Z}^D$ for a normal model with unknown mean and known variance $\sigma^2$.  We use Eq.~\ref{eq:statpartition} to treat the problem in terms of sufficient statitics. The statistic partition function breaks up by dimension: For each dimension with the flat improper prior:
\begin{equation} 
\varpi(\mu) = \sum_m \delta(\mu - m), 
\end{equation}
and sufficient statistic $t = N^{-1} \sum_i^N X_i$ the statistic partition function becomes the sum over  discrete prior values:
\begin{align}
z(t;N,\varpi ) &= \sum^{\infty}_{m=-\infty} q (t | m) = \left( \frac{N}{2 \pi} \right)^{1/2} \sum_{m=-\infty}^{\infty} e^{-N({t-m)^2}}\\
&=  \vartheta \left( t; r = e^{-\frac{2\pi^2}{N} } \right) = \sum_{m=-\infty}^{\infty} r^{m^2} e^{2 \pi t}
\intertext{Where $\vartheta$ is the Jacobi theta function with nome $r$. We can use the Jacobi triple product formula to write down an epxression for the log partition function:}
\log z(t;N,\varpi ) &= \sum_{m=1}^{\infty} \log(1- r^{2m}) + \log\left( 1+ r^{2m-1} e^{-i 2 \pi t} \right) + \log \left(1+ r^{2m-1} e^{i 2 \pi t}\right).\label{eqn:tripprod}
\end{align}
Assume (without loss of generality) that $m_0 =0$, then, because $|r| < 1$, we can Taylor expand the logarithm: 
\begin{align}
\expect{t|m_0}  \sum_{m=1}^{\infty} \log\left( 1+ r^{2m-1} e^{-i 2 \pi t} \right)  &= -  \sum_{m=1}^{\infty} \sum_{k=1}^{\infty} \frac{(-1)^k}{k} r^{(2m-1)k} \expect{t|m_0}   e^{-i 2 \pi k t},
\intertext{then compute the expectation term by term:}
&=-  \sum_{m=1}^{\infty} \sum_{k=1}^{\infty} \frac{(-1)^k}{k} r^{(2m-1)k+k^2},\\
 &=   \sum_{k=1}^{\infty} \frac{(-1)^k}{k} \frac{r^{k^2+k}}{r^{2 k}-1}.
\end{align}
This series is convergent, but converges slowly for very large $N$. We therefore must also develop a series for when $N\gg1$. We can use the Poisson resummation formula to convert the partition function into a sum over reciprocal space. First we have to subtract off the singularity at zero, by adding a piece to the summand that can be explicitly summed. Then we can extend this function to both positive and negative integers:
\begin{align}
\expect{t|m_0}  \sum_{m=1}^{\infty} \log \left( 1+ r^{2m-1} e^{-i 2 \pi t} \right)  &= \frac{1}{2} \sum_{k=1}^{\infty}  \frac{\cos(\pi k)}{k^2} \frac{k}{\sinh  \left( \frac{2 \pi^2 k}{N} \right)} e^{-\frac{2 \pi^2 k^2}{N}} \\
&= \frac{N}{48} + \frac{N}{4 \pi^2} \sum_{k=1}^{\infty} 	\frac{\cos(\pi k)}{k^2}
	\left( \frac{2 \pi^2 k}{n \sinh  \left( \frac{2 \pi^2 k}{N} \right)}		 e^{-\frac{2 \pi^2 k^2}{N}} - 1 \right)\\
&= \frac{N}{48} + \frac{1}{4} - \frac{\pi ^2}{12 N}  + \frac{N}{8 \pi^2} \sum_{k=-\infty}^{\infty} 	\frac{\cos(\pi k)}{k^2} 
	\left( \frac{2 \pi^2 k}{n \sinh  \left( \frac{2 \pi^2 k}{N} \right)}		 e^{-\frac{2 \pi^2 k^2}{N}} - 1 \right).
 \end{align}
The sum can now be represented as the sum of the Fourier transform of the summand. 
At large $N$, even the first term is exponentially small, and the whole sum can be ignored, leaving:
\begin{align}
\expect{t|m_0}  \sum_{m=1}^{\infty} \log \left( 1+ r^{2m-1} e^{-i 2 \pi t} \right)  &= 
\begin{cases} 
      \sum_{k=1}^{\infty} \frac{(-1)^k}{k} \frac{r^{k^2+k}}{r^{2 k}-1} & \text{for all\ } N\\
      \frac{N}{48} + \frac{1}{4} - \frac{\pi ^2}{12 N} & N > 10^2 \\
   \end{cases}
 \end{align}
Similarly we have for the Euler function piece of the triple product
\begin{align}
\sum_{m=1}^{\infty} \log(1- r^{2m}) =
\begin{cases}
\sum_{k=1}^{\infty}\frac{1}{k} \frac{r^{2 k}}{1-r^{2 k}} & \text{for all\ } N\\
-\frac{N}{24}+\frac{\pi ^2}{6 N}+\frac{1}{2} \log \left(\frac{N}{2 \pi }\right) & N > 50.
\end{cases}
\end{align}
%
The other term in Eq.~\ref{eqn:tripprod} can be computed in the same way. The free energy can then be computed using Eq.~\ref{eq:statpartition} and the computation of other derived quantities (learning capacity, GPI prior \textit{etc}) is left to the reader. 

%

\subsubsection{Normal model unknown mean and variance}
The  likelihood for the normal model is defined by Eq.~\ref{eqn:normalmodel} with parameters ${\bm \theta}=(\vec{\mu},\sigma)$ with support $\vec{\mu}\in \mathbb{R}^D$ and $\sigma\in \mathbb{R}_{>0}$.  The cross entropy is:
\begin{equation}
H({\bm \theta};{\bm \theta}_0) \equiv \textstyle\frac{D}{2} \left[\log 2\pi \sigma^2 + \textstyle\frac{\sigma^2_0}{\sigma^2}\right]+\textstyle\frac{1}{2\sigma^2} (\vec{\mu}-\vec{\mu}_0)^2, 
\end{equation}
where the true distribution is $X\sim q(\vec{x}|{\bm \theta}_0)$ and the determinant of the Fisher information matrix is:
\begin{equation}
\det {\bm I} = 2\sigma^{-2(D+1)}.
\end{equation}
The scaled Jeffreys prior (Eq.~\ref{eqn:scaledJeff}) is therefore: 
\begin{equation}
\rho = \sqrt{2}(\textstyle\frac{N}{2\pi \sigma^2})^{(D+1)/2}.
\end{equation}
We will assume $w$ matches the asymptotic form:
\begin{equation}
w = c\, \sigma^{-D-1}.
\end{equation}
Note that $w$ must have units of inverse length to the $D+1$ power in order to give the evidence the correct units. Due to translation symmetry in $\mu$, $w$ must be a function of $\sigma$ only.
The partition function is
\begin{equation}
\log Z(X^N) \sim \log c - \textstyle\frac{DN}{2} \log 2\pi \sigma^2 + \textstyle\frac{D}{2}\log \textstyle\frac{2\pi}{N}-\log 2+\log \Gamma(\textstyle \frac{DN}{2}) -\textstyle\frac{DN}{2}\log \textstyle\frac{\chi^2_{D(N-1)}}{2}
\end{equation}
where $\chi^2_{D(N-1)}$ is a $D(N-1)$-dimensional chi-squared random variable. The expected log partition function is:
\begin{equation}
\overline{\log Z}({\bm \theta};N) =  -NH_0({\bm \theta}) + \log c + \textstyle\frac{DN}{2}  + \textstyle\frac{D}{2}\log \textstyle\frac{2\pi}{N}-\log 2+\log \Gamma(\textstyle \frac{DN}{2}) -\textstyle\frac{DN}{2}\psi(\textstyle\frac{D(N-1)}{2})
\end{equation}
where $H_0$ is the true entropy and  $\psi$ is the polygamma function. The learning capacity is:
\begin{equation}
\overline{C}({\bm \theta};N,w) = \textstyle\frac{D}{2} + N^2(\textstyle\frac{D}{2})^2\psi^{(1)}(\textstyle\frac{DN}{2}) - 2N^2(\textstyle\frac{D}{2})^2\psi^{(1)}(\textstyle\frac{D(N-1)}{2})- N^3(\textstyle\frac{D}{2})^3\psi^{(2)}(\textstyle\frac{D(N-1)}{2}),
\end{equation}
where $D$ is both the dimension of mean parameter and the model. The
unknown normalization of $w$ is: 
\begin{equation}
\log c = \textstyle\frac{D+1}{2}\log \textstyle\frac{N}{2\pi} + \textstyle\frac{1}{2}\log 2 - {\cal K}
\end{equation}
written in terms of the effective complexity:
\begin{equation}
\inliner{{\cal K} =  \textstyle\frac{1}{2}\log \textstyle\frac{N}{2\pi}-\textstyle\frac{1}{2}\log 2 - \textstyle\frac{DN}{2}\log \textstyle\frac{N}{N+1} - N\log \Gamma[ \textstyle\frac{D(N+1)}{2} ] + (N+1)\log \Gamma[ \textstyle\frac{DN}{2} ]+  \textstyle\frac{D(N+1)N}{2} \left[ \psi(\textstyle\frac{DN}{2} -\psi(\textstyle\frac{D(N-1)}{2}))\right]}, \label{eqn:nmumuv}
\end{equation}
which is used to define the GPI prior $w$ using Eq.~\ref{eqn:calK}.

\subsection{Uniform distribution}
\label{sec:appuni}

In the exponential mixture model example, the non-regular model showed reduced learning capacity at the singularity but non-regular models can also have increased learning capacity as well. To illustrate this phenomenon, consider a continuous version of the German Tank problem, estimation of the end point of a uniform distribution \cite{goodman1954some}.

The  likelihood for the normal model is defined:
\begin{equation}
q(x|{\bm \theta}) \equiv \begin{cases} L^{-1}, &0\le x \le L\\
0, & \text{otherwise} \end{cases}, 
\end{equation}
where the parameter ${\bm \theta} \equiv (L)$ with support $L\in \mathbb{R}_{>0}$.
The cross entropy is:
\begin{equation}
H({\bm \theta};{\bm \theta}_0) = \begin{cases} \log L, & L_0 \le L\\
\infty, & \text{otherwise} \end{cases},  
\end{equation}
which is minimized at $L=L_0$ but neither the first nor second derivative is defined at this point and therefore the Fisher information matrix cannot be defined. We can still infer the dependence of the $w$ by symmetry and dimensional analysis:
\begin{equation}
w = c\, L^{-1}. \label{eqn:uw}
\end{equation}
The partition function is
\begin{equation}
\log Z(X^N) \sim \log c-N \log L-\log N-N \log Y,
\end{equation}
where $Y$ is the maximum of $N$ uniformly-distributed random variables on the interval $[0,1]$. The CDF for $Y$ is the $N$th power of the CDF for a single uniformly-distributed random variable. 
The expected log partition function is:
\begin{equation}
\overline{\log Z}({\bm \theta};N) = -NH_0({\bm \theta}) +\log c - \log N + 1.
\end{equation}
The learning capacity is:
\begin{equation}
\overline{C}({\bm \theta};N) = 1.
\end{equation}
The
unknown normalization of $w$ is: 
\begin{equation}
\log c = \log N - N\log (1+N^{-1}) -1,
\end{equation}
which can be plugged into Eq.~\ref{eqn:uw} to calculation the GPI prior $w$.

\subsection{Gamma model}

\label{app:Exp}

The  likelihood for the Gamma model is defined:
\begin{equation}
q(x|{\bm \theta}) \equiv \textstyle\frac{\beta^\alpha}{\Gamma(\alpha)} x^{\alpha-1}e^{-\beta x}
\end{equation}
where the parametersare ${\bm \theta} \equiv (\beta)$ with support $\beta \in \mathbb{R}_{>0}$. Note that the expontial model corresponds to $\alpha = 1$, the Uniform model corresponds to $\alpha \rightarrow 0$ and the normal model with unknown variance known mean corresponds to $\alpha = \textstyle\frac{1}{2}$, after the transformation $x\rightarrow x^{1/\alpha}$. 

The cross entropy is:
\begin{equation}
H({\bm \theta};{\bm \theta}_0) = -\alpha \log \beta + \log \Gamma(\alpha) - (\alpha-1)\left[\psi(\alpha)-\log \beta \right] + \textstyle\frac{\beta \alpha}{\beta_0}, 
\end{equation}
where the true distribution is $X\sim q(x|{\bm \theta}_0)$ and the determinant of the Fisher information matrix is:
\begin{equation}
\det {\bm I} = \frac{\alpha}{\beta^2}.
\end{equation}
The scaled Jeffreys prior (Eq.~\ref{eqn:scaledJeff}) is therefore: 
\begin{equation}
\rho = (\textstyle\frac{\alpha N }{2\pi\beta^2})^{1/2}.
\end{equation}
We will assume $w$ matches the asymptotic form:
\begin{equation}
w = c\, \beta^{-1}.
\end{equation}
The partition function is
\begin{equation}
\log Z(X^N) \sim -\alpha N\log Y^{\alpha N}+(\alpha-1) \sum_{i=1}^N\log Y^\alpha_i-N\log \beta_0 + \log \Gamma(\alpha N)-N\log \Gamma(\alpha),
\end{equation}
where $Y^{I}$ is a Gamma-distributed random variable with unit scale and shape $I$. 
The expected log partition function is:
\begin{equation}
\overline{\log Z}({\bm \theta};N) =  -NH_0({\bm \theta}) + \log \Gamma(N \alpha )-\alpha N \psi(N\alpha) + N\alpha,
\end{equation}
where $H_0$ is the true entropy and  $\psi$ is the polygamma function. The learning capacity is:
\begin{equation}
\overline{C}({\bm \theta};N) = -(\alpha N)^2[\psi^{(1)}(\alpha N)+\alpha N\psi^{(2)}(\alpha N)],
\end{equation}
where  $\psi$ is the polygamma function.
The
unknown normalization of $w$ is: 
\begin{equation}
\log c = \textstyle \frac{1}{2} \log \textstyle \frac{\alpha N}{2\pi}  - {\cal K}
\end{equation}
which can be re-written as an effective complexity:
\begin{equation}
{\cal K} = \textstyle \frac{1}{2} \log \textstyle \frac{\alpha N}{2\pi}+ (1+N) \log\Gamma(\alpha N) + -  N \log \Gamma(\alpha (1 + N)) - \alpha N (1 + N) (\psi[ \alpha N] - \psi[\alpha (1 + N)] ),\label{eqn:emum}
\end{equation}
to define the GPI prior $w$ using Eq.~\ref{eqn:calK}.

\subsection{GPI prior for discrete parameter manifolds}
\label{sec:discretew}

For discrete parameter manifolds two competing and well-established methods exist for choosing a prior:
(i) A literal interpretation of the principle of indifference would seem to imply that all parameter values are given equal weight. (ii) Alternatively, we can consider the continuous parameter limit where the prior can be chosen to give consistent results with the Jeffreys prior. Both approaches have desirable properties in different  contexts \cite{Berger2012}. GPI provides an elegant resolution to this conflict:
 When the discrete nature of the parameter manifold can be statistically resolved, the GPI prior prior assigns equal weight to all discrete parameter values (i), whereas, if the discreteness of the space cannot be statistically resolved, the large $N$ limit gives rise to a Jeffreys prior (ii):\begin{equation}
w = \begin{cases}
 \rho\ e^{-K} \prod_i \Delta \theta^i, & \Delta \theta^i \ll \delta \theta^i \\ 
1, & \Delta \theta^i \gg \delta \theta^i
\end{cases}\label{eqn:discretew}
\end{equation}
where $\Delta \theta^i$ is the lattice spacing and the statistical resolution $\delta \theta^i$ is defined in Eq.~\ref{eqn:resolution}. The GPI prior for a normal model with a discrete mean can be computed exactly and is described in the appendix (\ref{App:discra}).

%
%
%
%
%
%
%

\bibstyle{PRX}

\end{document}